\documentclass[10pt,a4paper]{article}
\usepackage{amsmath}
\usepackage{amssymb}
\usepackage{epsf}
\usepackage{amscd}
\usepackage[all]{xy}
\usepackage{enumerate}

\input xy
\xyoption{all}

\numberwithin{equation}{section}

\newtheorem{theo}{Theorem}[section]

\newtheorem{prop}[theo]{Proposition}

\newtheorem{cor}[theo]{Corollary}

\newtheorem{lemma}[theo]{Lemma}

\newtheorem{dfn}[theo]{Definition}

\newtheorem{defprop}[theo]{Definition-Proposition}

\newtheorem{rquee}[theo]{{\it Remark}}
\newenvironment{rmk}{\begin{rquee} \normalfont}{\end{rquee}}

\newtheorem{csqcee}[theo]{{\it Consequence}}

\newtheorem{exple}[theo]{{\it Example}}
\newenvironment{ex}{\begin{exple} \normalfont}{\end{exple}}

\newtheorem{pty}[theo]{Property}

\newtheorem{notationse}[theo]{{\it Notation}}
\newenvironment{nota}{\begin{notationse} \normalfont}{\end{notationse}}

\newenvironment{proof}{{\flushleft{\it Proof: }} \rm}{\Qed \\}

\newenvironment{proofproposition}{{\flushleft{\it Proof of Proposition}} \rm}{\Qed \\}

\newcommand{\bd}{\begin{dfn}}
\newcommand{\ed}{\end{dfn}}
\newcommand{\bp}{\begin{prop}}
\newcommand{\edp}{\end{defprop}}
\newcommand{\bdp}{\begin{defprop}}
\newcommand{\ep}{\end{prop}}
\newcommand{\bt}{\begin{theo}}
\newcommand{\et}{\end{theo}}
\newcommand{\bc}{\begin{cor}}
\newcommand{\ec}{\end{cor}}
\newcommand{\be}{\begin{ex}}
\newcommand{\ee}{\end{ex}}
\newcommand{\bexo}{\begin{exo}}
\newcommand{\eexo}{\end{exo}}
\newcommand{\bl}{\begin{lemma}}
\newcommand{\el}{\end{lemma}}
\newcommand{\br}{\begin{rmk}}
\newcommand{\er}{\end{rmk}}
\newcommand{\bpy}{\begin{pty}}
\newcommand{\epy}{\end{pty}}
\newcommand{\bpf}{\begin{proof}}
\newcommand{\epf}{\end{proof}}
\newcommand{\bpfl}{\begin{proof}}
\newcommand{\epfl}{\end{proof}}
\newcommand{\bpfp}{\begin{proofproposition}}
\newcommand{\epfp}{\end{proofproposition}}
\newcommand{\bn}{\begin{nota}}
\newcommand{\en}{\end{nota}}

\newcommand{\qed}{\hfill \mbox{$\square$}}
\newcommand{\Qed}{\hfill \mbox{$\blacksquare$}}

\newcommand{\ppq}{\leqslant}
\newcommand{\pgq}{\geqslant}

\newcommand{\m}{\frak{m}}

\newcommand{\cc}{\mathbb{C}}
\newcommand{\zz}{\mathbb{Z}}

\newcommand{\Bb}{\mathbb{B}}

\newcommand{\rad}{\mathrm{rad}}
\newcommand{\soc}{\mathrm{soc}}
\newcommand{\Top}{\mathrm{top}}
\newcommand{\id}{\mathrm{id}}
\newcommand{\im}{\mathrm{Im}}
\renewcommand{\ker}{\mathrm{Ker}}
\newcommand{\tr}{\mathrm{tr}}
\newcommand{\qtr}{\mathrm{qtr}}
\newcommand{\qdim}{\mathrm{qdim}}

\newcommand{\spn}{\mathrm{span}}
\newcommand{\car}{\mathrm{char}}
\newcommand{\gr}{\mathrm{deg}_{\ot}}
\newcommand{\Hom}{\mathrm{Hom}}

\newcommand{\End}{\mathrm{End}}

\newcommand{\Ext}{\mathrm{Ext}}
\newcommand{\A}{\mathcal{A}}

\newcommand{\C}{\mathcal{C}}
\newcommand{\D}{\mathcal{D}}

\newcommand{\ot}{\otimes}

\newcommand{\nequiv}{\equiv\hspace*{-9pt}/\ }
\newcommand{\proj}{2j+u-1\equiv 0 \pmod{d}}
\newcommand{\nproj}{2j+u-1\nequiv 0 \pmod{d}}
\newcommand{\vv}[1]{\tilde{H}_{#1}}
\newcommand{\elt}[3]{\gamma_{#2}^{#3}\tilde{H}_{#1,#2}}

\def\land{\Lambda _{n,d}}

\def\dland{\mathcal{D}(\Lambda _{n,d})}
\def\m{\frak{m}}

\newcommand{\set}[1]{\{ #1 \}}

\newcommand{\mat}[2]{\mathcal{M}_{#1}(#2)}
\newcommand{\rep}[1]{\langle{#1}\rangle}

\newcommand{\abs}[1]{|{#1}|}

\newcommand{\mx}[1]{\begin{pmatrix}#1\end{pmatrix}}

\title{\textsc {Representation theory of the Drinfel'd doubles of a family of Hopf algebras}}
\author{K. Erdmann \and E.L. Green\thanks{Partially supported by an
    NSA grant} \and N. Snashall \and R. Taillefer\thanks{Partially
    supported by an LMS  Scheme 4 grant and a Lavoisier grant from the
    French Ministry of Foreign Affairs}}
\date{}

\begin{document}

\maketitle

\begin{abstract}We investigate the Drinfel'd doubles $\dland$ of a certain family of Hopf algebras. We determine their simple modules and their indecomposable projective modules, and 
 we obtain a presentation by quiver and relations of these Drinfel'd doubles, from which we deduce properties of their representations, including the Auslander-Reiten quivers of the $\dland$. We then determine decompositions of the tensor products of most of the representations described, and in particular give a complete description of the tensor product of two simple modules. This study also leads to explicit examples of Hopf bimodules over the original Hopf algebras.
\end{abstract}

{\flushleft\bf Mathematics Subject Classification (2000):} 17B37, 06B15, 81R50, 16W30, 16W35, 16G20, 16G70, 18D10.


\section{Introduction}

We study the representation theory of some Drinfel'd
doubles. The Drinfel'd double of a finite-dimensional Hopf algebra was
defined by Drinfel'd in order to provide solutions to the quantum
Yang-Baxter equation arising from statistical mechanics. 

Representations of a Hopf algebra (up to isomorphism) form a ring in
which the product is given by the tensor product over the base field, and
in the case of the Drinfel'd double (and any quasitriangular Hopf
algebra) this ring is commutative. 

However, not very much is known about the representations of the
Drinfel'd double of a nonsemisimple Hopf algebra in general. S. Witherspoon studied the Drinfel'd
double of the algebra of a finite group in positive characteristic
\cite{witherspoon}. She proved in particular that the Green ring of
the Drinfel'd double of a group algebra decomposes as a product of
ideals associated to some subgroups of the original
group. H-X. Chen 
gives a complete list of simple modules over the Drinfel'd doubles of
the Taft algebras in \cite{chen2}.
D.E. Radford in \cite{Rd2} studies the simple modules over the
Drinfel'd double of a finite-dimensional Hopf algebra, and
characterises them in some cases (which include the Taft
algebras and quiver Hopf  algebras); he then uses them to construct twist oriented quantum
algebras, which give rise to invariants of oriented knots and
links. Moreover, he establishes that all modules over the Drinfel'd
double can be viewed as submodules of some particular Yetter-Drinfel'd
modules (constructed from modules over the original Hopf algebra).

An important property of the algebra underlying the Drinfel'd double of a finite-dimensional
 Hopf algebra $H$ can be deduced from work of D.E. Radford and
 R. Farnsteiner (\cite[Corollary 2 and Theorem 4]{Rd} and
 \cite[Proposition 2.3]{farnsteiner}): the Drinfel'd double $\mathcal{D}(H)$
is a symmetric algebra.\label{symmetric}

In this paper, we study in detail the representation theory of the
Drinfel'd doubles of the algebras $\land,$ that is, of the duals of the
extended Taft algebras (defined in the first section). The algebras $\land$
are interesting, since they are finite-dimensional Hopf
algebras which are neither commutative nor cocommutative, and which
are not quasitriangular, with
antipodes of arbitrarily high order ($2n$). However, they are monomial (that is, all the relations are given by paths)
and have finite
representation type, and as such are
therefore very suitable to study.

In the first section, we study the representations of these Drinfel'd doubles,
first describing the projective modules and the simple modules (see \ref{simple}, \ref{projective}, \ref{simples} and \ref{projectivesummary}). We give a complete presentation by quiver and relations, and this shows that the algebras are of tame representation type (see Section \ref{sectionclassification}; see also \cite[4.4]{benson} for definitions of representation types).
These quivers are similar to some quivers which occur in the
representation theory of blocks of reduced enveloping algebras in
characteristic $p$ (see \cite{FS}). However, the situation in our case is quite different: as
we see later on in the paper, there exist representations which are
periodic of period greater than 2 in general (which is not the case
for reduced enveloping algebras). Moreover, the conditions on the
characteristic of the field  in our context and in the
case of reduced enveloping algebras are distinct; in the case of reduced
enveloping algebras, the cycles which occur in the quiver have a
length which is a power of the characteristic, whereas in our case the
characteristic does not divide the length of the cycles but is otherwise arbitrary.

In the following section we describe all the indecomposable
representations of $\dland,$ and the Auslander-Reiten quiver of these
algebras. We also determine which of these representations are
splitting trace modules and more precisely give the quantum dimensions
of the modules (see Section \ref{sectionsplittingtracemodules}, where these concepts are also defined).  We then study tensor products of representations. We
first give a complete and explicit decomposition of the tensor product of two
simple modules, and in particular prove that this is always semisimple modulo some projective direct summands (see \ref{tensorsimples}). We  then describe the tensor products of many other modules up to
projectives. This gives a large part of the structure of the Green
ring of $\dland$ modulo projectives.  Finally, we give some examples of Hopf bimodules over the
original Hopf algebras $\land$, obtained via the equivalence of
categories between Hopf bimodules over $\land$ and modules over
$\dland.$
We conclude with an appendix in which we outline the proof of the classification of
the representations of $\dland.$

This work opens many questions on the representation theory of
Drinfel'd doubles, in particular the possibility of finding some
general properties of tensor products of simple modules over the Drinfel'd
doubles of some Hopf algebras, using D.E. Radford's results \cite{Rd2}.

For general facts about representations of algebras, we refer to \cite{ARS}, \cite{benson} and \cite{erdmann}. We will refer to those more precisely in some parts of the paper.

\ 

{\sc Acknowledgements:} The second and third authors thank the University of Oxford and the last author thanks the Virginia Polytechnic
Institute and State University for their
kind hospitality.


\section{Projectives and quiver for $\dland$}

\subsection{Preliminaries}

Throughout the paper, $k$ is  an algebraically closed field. 

The algebra $\land$ is described by quiver and relations;  we refer to \cite[III.1]{ARS} and \cite[4.1]{benson} for definitions and properties relating to quivers and relations. The quiver
is cyclic,
$$
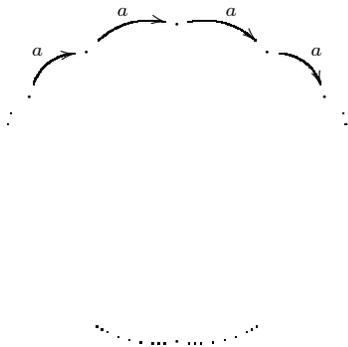
\xymatrix@=.01cm{
&&&&&&&&&&&\cdot\ar@/^.5pc/[rrrrrd]^{a}\\
&&&&&&\cdot\ar@/^.5pc/[rrrrru]^{a}&&&&&&&&&&\cdot\ar@/^.5pc/[rrrdd]^{a}\\\\
&&&\cdot\ar@/^.5pc/[rrruu]^{a}\ar@{.}@/_.3pc/[ldd]&&&&&&&&&&&&&&&&\cdot\ar@{.}@/^.3pc/[rdd]\\\\
&&&&&&&&&&&&&&&&&&&&&&&\\
\\
\\
\\\\\\\\\\\\\\\\
&&&&&&&&&&&&&&&&\\
&&&&&&&&&&&\cdot\ar@/_.3pc/@{.}[rrrrru]\ar@/^.3pc/@{.}[lllllu]
}$$
 with $n$ vertices $e_0,\ldots, e_{n-1}$ and $n$ arrows $a_0,\ldots,a_{n-1}$,
where the arrow $a_i$ goes from the vertex $e_i$ to the vertex
$e_{i+1},$ and we factor by the ideal generated by all paths of length
$d\pgq 2$. We shall denote by $\gamma_i^m$ the path $a_{i+m-1}\ldots
a_{i+1}a_i$ (read from right to left), that is, the path of length $m$ starting at
the vertex $e_i$. 

When $d$ divides $n,$ this algebra is a Hopf algebra, and in fact the
condition $d\mid n$ is a necessary and sufficient condition for
$\land$ to be a  Hopf algebra when $\car k=0$ ({see}
\cite{cibils,chen1}). 
This Hopf algebra can actually be considered over more general fields,
and in this paper, we assume only that the characteristic of $k$ does
not divide $n$ (to ensure existence of roots of unity).

 We fix  a primitive 
 $d^{th}$ root of unity $q$ in $k.$  The formulae
$$
\begin{array}{lclcl}
\varepsilon(e_i)=\delta_{i0}& &\Delta(e_i)=\sum_{j+\ell=i}e_j\ot e_\ell&&S(e_i)=e_{-i}\\
\varepsilon(a_i)=0 & & \Delta(a_i)=\sum_{j+\ell=i}(e_j\ot
a_\ell+q^\ell a_j\ot
e_\ell)&&
S(a_i)=-q^{i+1}a_{-i-1}
\end{array}
$$ determine the Hopf algebra structure of $\land.$ 

We want to study the Drinfel'd double of $\land.$ Recall:

\bd[{see} \cite{kassel,montgomery}]\label{doubledef}  Let $H$ be a Hopf
algebra. The {\bf Drinfel'd double} of $H$ is the Hopf algebra which is
equal to $H^{*cop}\ot H$ as a coalgebra (ordinary tensor product of coalgebras), and whose product is defined
by $$(\alpha\ot h)(\beta\ot g)=\alpha \beta (S^{-1}h^{(3)}?h^{(1)})\ot
h^{(2)}g,$$ where $\beta (S^{-1}h^{(3)}?h^{(1)})$ is the map which
sends $x\in H$ to $\beta (S^{-1}h^{(3)}xh^{(1)})\in k$ and where we
have used the Sweedler notation $\Delta(h)=h^{(1)}\ot h^{(2)}$ for the
comultiplication.\ed

Therefore, we need to understand the dual $\land^{*cop}$; it is isomorphic as an algebra
to the {extended Taft algebra}\label{taft} $$\rep{G,X\mid
  G^n=1,X^d=0,GX=q^{-1}XG}.$$  (C. Cibils explains this for $n=d$ when
$\land$ is a selfdual Hopf algebra in \cite{cibils}, and
the general case is similar; if $\set{\check{\gamma}_i^m\,\mid\,
  i\in\zz_n, 0\ppq m\ppq d-1}$ denotes the dual basis of
$\land^{*cop},$ the correspondence is determined by
$\check{e}_i\mapsto G^i$ and
$\check{a}_i\mapsto G^iX$). Its Hopf algebra structure is defined
by $$
\begin{array}{lclcl}
\varepsilon(G)=1&&\Delta(G)=G\ot G&&S(G)=G^{-1}\\
\varepsilon(X)=0&&\Delta(X)=X\ot G+1\ot X &&S(X)=-XG^{-1}=-q^{-1}G^{-1}X.
\end{array}$$

\bn We shall need the following notation to describe some elements in
$\dland:$ for non-negative integers $m$ and $u$, we define the
{\bf $q$-integers} $(0)_q=0$ and $(m)_q=1+q+q^2+\ldots+q^{m-1},$
the {\bf $q$-factorials} $0!_q=1$ and $m!_q=(m)_q(m-1)_q\ldots (2)_q(1)_q,$  and
the {\bf $q$-binomial coefficients} $\begin{pmatrix}m\\u
\end{pmatrix}_{\!q}=\displaystyle{\frac{m!_q}{(m-u)!_qu!_q}}$
({see} for instance \cite{cibils,Rd3}, where some properties are also given).
\en

\br Using the notation above, we have $\Delta(\gamma_i^m)=\displaystyle{\sum_{
{\tiny\begin{array}{c}j\in \zz_n\\0\ppq v\ppq m
\end{array}}}
}\begin{pmatrix}m\\v \end{pmatrix}_{\!q}   q^{vj}\gamma_{i-j}^v\ot
\gamma_j^{m-v}$ and $\Delta(X^m)=\sum_{v=0}^m \mx{m\\v}_{\!q}X^v\ot G^vX^{m-v}.$\er

We can now describe the Drinfel'd double $\dland:$

\bp\label{doublepresentation} The Drinfel'd double $\dland$ is described as follows: as a
coalgebra, it is $\land^{*cop}\ot \land$. We write the basis elements
$G^iX^j\gamma_\ell^m,$ with $i,\ell \in \zz_n$ and $0\ppq j,m \ppq d-1$
(\emph{i.e.} we do not write the tensor product symbol). The following
relations determine the algebra structure completely: 

\begin{eqnarray*}
&&  G^n=1, X^d=0, \mbox{ and } GX=q^{-1}XG, 
\\ && \mbox{the product of elements $\gamma_\ell^m$ is the
usual product of paths},
 \\ && \gamma_\ell^mG=q^{-m}  G\gamma_\ell^m, \mbox{and}
\\&&  \gamma_\ell^mX= q^{-m} X\gamma_{\ell+1}^m -q^{-m}(m)_{\!q}
 \gamma_{\ell+1}^{m-1}+q^{\ell+1-m} (m)_{\!q}G\gamma_{\ell+1}^{m-1}.
\end{eqnarray*}

 \ep 

\bpf
The proof is straightforward, using the definition of the Drinfel'd
double (Definition \ref{doubledef}) and the structure of $\land^{*cop}.$
\epf

\bn\label{reps} Since some indices are described modulo $d$ and others modulo $n,$
we need to make a distinction. If $j$ is an element in $\zz_n,$ we
shall denote its representative modulo $d$ in $\set{1,\ldots,d}$ by
$\rep{j}$ and its representative modulo $d$ in $\set{0,\ldots, d-1}$ by
$\rep{j}^-.$\en


\subsection{Quiver of $\dland$}

The aim of this section is to describe the quiver of $\dland.$ To do
this, we first decompose $\dland$ into a (non-minimal) product of
algebras $\Gamma_0,$ $\ldots,$ $\Gamma_{n-1},$ 
and study each  of these algebras. We give a basis for $\Gamma_u$ for $u=0,\ldots,n-1,$ and
describe some indecomposable $\Gamma_u$-modules as a first step in the
construction of the indecomposable projective modules.
 This follows R. Suter's method in \cite{suter}
where he studies the representations of a finite-dimensional quotient
of $U_q(\frak{sl}_2(\cc))$ ({see} also \cite{xiao,patra}).


\subsubsection{First decomposition of $\dland$}

\bp The elements $E_u:=\frac{1}{n}\sum_{i,j\in
  \zz_n}q^{-i(u+j)}G^ie_j$, for $u \in \zz_n,$ are central orthogonal
idempotents, and $\sum_{u \in \zz_n}E_u=1.$ 

Therefore $\dland\cong\prod_{u\in\zz_n} \Gamma_u$ where $\Gamma_u=\dland E_u.$\ep

We shall now study $\Gamma_u.$


\subsubsection{Construction of  modules over $\dland$}

We now define some idempotents inside $\Gamma_u,$ which are not
central, but we will use them to describe a basis for $\Gamma_u.$

\bp Set $E_{u,j}=\sum_{v=0}^{\frac{n}{d}-1}e_{j+vd}E_u$, for $j\in\zz_d$. Then
$E_{u,j}E_{u,\ell}=\delta_{j\ell}E_{u,j}$ and $\sum_{j=0}^{d-1}E_{u,j}=E_u.$ 
We also have $E_{u,j}=E_{u,j'}$ iff $j\equiv j' \pmod{d}.$

Moreover, the following relations hold within $\Gamma_u$:
$$
GE_{u,j}=q^{u+j}E_{u,j}=E_{u,j}G \hspace{1cm}  XE_{u,j}=E_{u,j-1}X$$ $$ 
\gamma_\ell ^mE_{u,j}=\begin{cases} E_{u,j+m}\gamma_\ell ^m &\mbox{ if
    $\ell \equiv j \pmod{d}$}\\0&\mbox{ otherwise.} \end{cases}$$
\ep

We can now describe a basis for $\Gamma_u$ and  a grading on
$\Gamma_u,$ as follows:
$\Gamma_u=\bigoplus_{s=-d+1}^{d-1}(\Gamma_u)_s$ with
$(\Gamma_u)_s=\spn\set{X^t\gamma_j^mE_{u,j}\,\mid\, j\in \zz_n, 0\ppq m,t
  \ppq d-1, m-t=s}.$ 

Multiplication on the left and on the right by elements in $\dland$ respect this grading, with multiplication by $X$ on either side
reducing the degree by 1, multiplication by an arrow on either
side increasing the degree by 1 and multiplication by $G$ on either
side leaving the degree unchanged. Moreover $(\Gamma_u)_s$ is a sum of
eigenspaces for $G:$  if $y E_{u,j}$ is an element in
$(\Gamma_u)_s,$ we have 
$G\cdot y E_{u,j}=q^{s+j+u} y E_{u,j}$ and $
y E_{u,j}\cdot G=q^{j+u} y E_{u,j}.$

\ 

Now set $F_{u,j}:=\gamma_j^{d-1}E_{u,j}$ for $j\in\zz_n.$ Then:

\bp \label{firstmodules} The module $\Gamma_uF_{u,j}$ has the following form:
$$\xymatrix@R=.5cm{{\scriptstyle F_{u,j}}\ar@<+.5ex>[d]^X \\ \ar@<+.5ex>[d]^X
  \ar@<+.5ex>^{a_{j+d-2}}[u] \\ \ar@<+.5ex>^{a_{j+d-3}}[u] \ar@{.}[d] \\ \ar@<+.5ex>[d]^X \\
 {\scriptstyle H_{u,j}} \ar@<+.5ex>[d]^X
\ar@<+.5ex>^{a_{}}[u]  \\ {\scriptstyle\tilde{F}_{u,j}}\ar@<+.5ex>[d]^X \\
\ar@<+.5ex>^{a_{}}[u]\ar@{.}[d] \\
\ar@<+.5ex>[d]^X\\ {\scriptstyle\tilde{H}}_{u,j}
\ar@<+.5ex>^{a_{j}}[u]}$$
where $H_{u,j}:=X^{\rep{2j+u-1}-1}F_{u,j},$
$\tilde{F}_{u,j}:=X^{\rep{2j+u-1}^-}F_{u,j}$ and
$\tilde{H}_{u,j}:=X^{d-1}F_{u,j}$ (recall the notations
in \ref{reps}). 

In this diagram, the arrows represent the actions of $X$ and of the
arrows in the original quiver up to a nonzero scalar; the basis vectors are eigenvectors for the action of $G$.

Note that when
$\proj,$ the single arrow does not occur, the module is simple, and we
have $H_{u,j}=\vv{u,j}=X^{d-1}F_{u,j},$ and $\tilde{F}_{u,j}=F_{u,j}.$

\ep

In order to prove this proposition, we require the following lemma:

 \bl\label{R1} The element $\gamma_{d+j-m-1}^bX^mF_{u,j}$
is equal to
$$q^{-\frac{b(2m-b+1)}{2}}\frac{(m)!_{\!q}}{(m-b)!_{\!q}}\prod_{t=1}^b
(q^{2j+u-1-(m-t+1)}-1)\; X^{m-b}F_{u,j}$$ if $b<m$ and is
0 otherwise.\el
\bpfl This is proved by induction on $b.$ When $b=1,$ we take the
arrow on the left across the $X$'s, using the relations in Proposition
\ref{doublepresentation}.

\begin{eqnarray*}
a_{j-m-1+d}X^m\gamma_j^{d-1}E_{u,j}&=&
\left(q^{-1}Xa_{j-m+d}-q^{-1}(1-q^{2(j-m+d)+u+2})e_{j-m+d}\right)X^{m-1}\gamma_j^{d-1}E_{u,j}\\
&=&q^{-1}\left(q^{-1}X^2a_{j-m+d+1}-q^{-1}(1-q^{2(j-m+d)+u+2})X\right)X^{m-2}\gamma_j^{d-1}E_{u,j}\\&&-q^{-1}(1-q^{2(j-m+d-1)+u+2})X^{m-1}\gamma_j^{d-1}E_{u,j}\\
&=& q^{-2} X^2a_{j-m+d+1}X^{m-2}\gamma_j^{d-1}E_{u,j}\\&& +
\left(-q^{-2}-q^{-1}+q^{2j-2m+u-1}+q^{2j-2m+u}\right)X^{m-1}\gamma_j^{d-1}E_{u,j}\\
&=&\ldots\\
&=& q^{-m}X^ma_{j+d-1}\gamma_j^{d-1}E_{u,j}\\&&+\sum_{p=1}^m (-q^{-p}+q^{2j-2m+u+p-2})X^{m-1}\gamma_j^{d-1}E_{u,j},
\end{eqnarray*}
 with the first term in the last identity equal to 0. \epfl  

\bpfp {\bf\ref{firstmodules}:} Here we apply Lemma \ref{R1} with $b=1$ to see that 
the element  $a_{d+j-m-1}X^mF_{u,j}$ is equal to 0 if
$m\equiv 2j+u-1 \pmod{d}$ and is a nonzero multiple of
$X^{m-1}F_{u,j}$ otherwise.
\epfp

\br The simple projective modules that we have  described in
Proposition \ref{firstmodules} when $\proj$, are
distributed as follows in the whole algebra $\dland:$

When $d$ is odd,  there are $\frac{n}{d}$ such simples
in each algebra $\Gamma_u$. 

When $d$ is even, there are $\frac{2n}{d}$ such simples in $\Gamma_u$ if
 $u$ is odd, and none in $\Gamma_u$ for $u$ even.
\er

\bd \label{sigma} We define permutations of the indices in $\zz_n$ by
$\sigma_u(j)=d+j-\rep{2j+u-1}.$ Note that the arrow going up from $H_{u,j}$ in the diagram in Proposition \ref{firstmodules} is $a_{\sigma_u(j)}.$\ed

\br\label{sigmapty} We can easily see that $\sigma_u(j)=j$ if and only if $\proj,$ and
that if $\nproj,$ then $\sigma_u^2(j)=j+d$ and so $\sigma_u$ has order
$\frac{2n}{d}.$\er

We shall now define larger modules (and we will see later that they
are a full set of representatives of the indecomposable projective
modules).

\subsubsection{Indecomposable projective modules and decomposition of $\dland$}

\bl If $\nproj ,$ then there exists an element $K_{u,j}$ homogeneous of
degree $d-\rep{2j+u-1}^--1$ such that
$H_{u,j}=a_{\sigma_u(j)-1}K_{u,j}.$\el

\bpfl Consider $H_{u,j}= X^{\rep{2j+u-1}-1}\gamma_j^{d-1}E_{u,j}$. We first take one arrow across the $X$'s; there exist scalars
$\alpha_1,\ldots, \alpha_{2j+u-2}$ in $k$ such that:
\begin{eqnarray*}
H_{u,j}&=& X^{\rep{2j+u-1}-1}\gamma_j^{d-1}E_{u,j}\\
&=& qX^{\rep{2j+u-1}-2}a_{j+d-1}X\gamma_j^{d-2}E_{u,j} + \alpha_1
X^{\rep{2j+u-1}-2}\gamma_j^{d-2}E_{u,j} \\
&=& q^2 X^{\rep{2j+u-1}-3}a_{j+d}X^2\gamma_j^{d-2}E_{u,j}
+(\alpha_1+\alpha_2)X^{\rep{2j+u-1}-2}\gamma_j^{d-2}E_{u,j} \\
&=&\ldots\\
&=& q^{2j+u-2 }
a_{\sigma_u(j)-1}X^{\rep{2j+u-1}-1}\gamma_j^{d-2}E_{u,j} + (\alpha_1+\ldots+\alpha_{2j+u-2})X^{\rep{2j+u-1}-2}\gamma_j^{d-2}E_{u,j}. 
\end{eqnarray*} We now repeat the process on the second term of the
last identity, and continue until there is an arrow in front of all
the terms; so there exist scalars $\beta'_i$, $\beta''_i$ and $\beta_i$ such that
the following identities hold:

\begin{eqnarray*}
H_{u,j}&=& q^{2j+u-2 }
a_{\sigma_u(j)-1}X^{\rep{2j+u-1}-1}\gamma_j^{d-2}E_{u,j} +
\beta_1'X^{\rep{2j+u-1}-2}\gamma_j^{d-2}E_{u,j} \\
&=& q^{2j+u-2 }
a_{\sigma_u(j)-1}X^{\rep{2j+u-1}-1}\gamma_j^{d-2}E_{u,j} \\&&+
\beta_1'
\left(q^{2j+u-3}a_{\sigma_u(j)-1}X^{\rep{2j+u-1}-2}\gamma_j^{d-3}E_{u,j}
+\beta_2''X^{\rep{2j+u-1}-3}\gamma_j^{d-3}E_{u,j}\right) \\
&=& q^{2j+u-2 }
a_{\sigma_u(j)-1}X^{\rep{2j+u-1}-1}\gamma_j^{d-2}E_{u,j} + \beta_1
a_{\sigma_u(j)-1}X^{\rep{2j+u-1}-2}\gamma_j^{d-3}E_{u,j}\\&&+
\beta_2'X^{\rep{2j+u-1}-3}\gamma_j^{d-3}E_{u,j}\\
&=& \ldots\\
&=& q^{2j+u-2 }
a_{\sigma_u(j)-1}X^{\rep{2j+u-1}-1}\gamma_j^{d-2}E_{u,j}
+a_{\sigma_u(j)-1}
\sum_{p=2}^{\rep{2j+u-1}}\beta_{p-1}X^{\rep{2j+u-1}-p}\gamma_j^{d-p-1}E_{u,j}\\
&=& a_{\sigma_u(j)-1} K_{u,j} 
\end{eqnarray*} with $K_{u,j}$ nonzero and homogeneous of degree $d-\rep{2j+u-1}-1=d-\rep{2j+u-1}^--1.$

 \epfl

\bd If $\proj,$ set $K_{u,j}=F_{u,j} (=\tilde{F}_{u,j}).$  Note that
it is homogeneous of degree $d-1.$\ed

Now consider $\Gamma_u K_{u,j}.$ The following result is immediate:  

\bdp\label{simple} Assume that $\proj.$ Define $L(u,j):=\Gamma_uK_{u,j};$
this module has the following structure: 
$$\xymatrix@R=.5cm{&\mbox{degree}&\mbox{$G$-eigenvalue}&\mbox{length}\\
{\scriptstyle F_{u,j}=K_{u,j}=\tilde{F}_{u,j}}\ar@<+.5ex>[d]^X &d-1&q^{j+u-1}&\ar@{--}[dd]\\ 
\ar@<+.5ex>[d]^X \ar@<+.5ex>^{a_{j+d-2}}[u] &\ar@{.}[dd]&\ar@{.}[dd]\\ 
\ar@<+.5ex>^{a_{j+d-3}}[u] \ar@{.}[d]&&&d-1\ar@{--}[dd] &&\\ 
\ar@<+.5ex>[d]^X &&&&\\
{\scriptstyle\tilde{H}}_{u,j=H_{u,j}}\ar@<+.5ex>^{a_{j}}[u]&0&q^{j+u}&&&}$$

\edp
 
In the case $\nproj,$ we obtain the following structure:

\bp\label{projective} Assume that $\nproj.$ The module $\Gamma_uK_{u,j}$ has the following structure:

$$\xymatrix@R=.5cm@=.5cm{&&&\mbox{degree}&\mbox{$G$-eigenvalue}&\mbox{length}\\
&  {\scriptstyle F_{u,j}}  \ar@<+.5ex>[d]^X \ar@{.}[rr] &&  {\scriptstyle d-1 } \ar@{.}[r] &  q^{j+u-1} \ar@{.}[r] &
\ar@{--}[d]    \\
&  \ar@<+.5ex>^{a_{j+d-2}}[u]  \ar@{.}[d]  & &\ar@{.}[d]&\ar@{.}[d]&   {\scriptstyle \rep{2j+u-1}-1} \ar@{--}[dd]     \\
&   \ar@<+.5ex>[d]^X  &&&  \\
&  {\scriptstyle H_{u,j}} \ar@<+.5ex>^{a_{\sigma_u(j)}}[u]
\ar@<+.5ex>[dr]^X \ar@{.}[rr]&& {\scriptstyle  d-\rep{2j+u-1}} \ar@{.}[r] &  q^{-j+1} \ar@{.}[r] &  \\
{\scriptstyle K_{u,j}}   \ar@<+.5ex>^{a_{\sigma_u(j)-1}}[ur]
  \ar@<+.5ex>[d]^X &&  {\scriptstyle \tilde{F}_{u,j}} \ar@<+.5ex>[d]^X\ar@{.}[r]
  & {\scriptstyle  d-\rep{2j+u-1}^--1}\ar@{.}[r]  &   q^{-j} \ar@{.}[r] &   \ar@{--}[d]   \\
\ar@<+.5ex>^{a_{\sigma_u(j)-2}}[u]\ar@<+.5ex>^{a_{\sigma_u(j)-2}}[urr]   \ar@{.}[d] &&
\ar@<+.5ex>^{a_{\sigma_u(j)-2}}[u] \ar@{.}[d]&\ar@{.}[d]&\ar@{.}[d]&{\scriptstyle d-\rep{2j+u-1}^--1}\ar@{--}[dd]  \\
\ar@<+.5ex>[d]^X      &&    \ar@<+.5ex>[d]^X      &&& \\
{\scriptstyle D_{u,j}}    \ar@<+.5ex>^{a_{j}}[u]   \ar@<+.5ex>^{a_{j}}[urr]  \ar@<+.5ex>[dr]^X   &    &{\scriptstyle\tilde{H}_{u,j}} \ar@<+.5ex>^{a_{j}}[u] \ar@{.}[r] &  {\scriptstyle  0}  \ar@{.}[r] &   q^{j+u} \ar@{.}[r]  & \\
&  {\scriptstyle \tilde{K}_{u,j}}  \ar@<+.5ex>[d]^X     \ar@<+.5ex>^{a_{j-1}}[ur] \ar@{.}[rr] && {\scriptstyle -1}\ar@{.}[r] &
q^{u+j-1}\ar@{.}[r] & \ar@{--}[dd]  \\
&\ar@<+.5ex>^{a_{j-2}}[u]   \ar@{.}[d] &&\ar@{.}[d]&\ar@{.}[d]  \\
&   \ar@<+.5ex>[d]^X   &&&&  {\scriptstyle \rep{2j+u-1}-1}  \ar@{--}[d]  \\
& {\scriptstyle \tilde{D}_{u,j}}  \ar@<+.5ex>^{a_{\sigma_u(j)-d}}[u]\ar@{.}[rr] && {\scriptstyle -\rep{2j+u-1}}\ar@{.}[r]  &  q^{-j+1} \ar@{.}[r]& 
}$$\ep 

To prove this, we require the following two lemmas:

\bl\label{commut} We have 
$\gamma_{\sigma_u(j)-s-1}^tX^sK_{u,j}=\sum_{b=s-t+1}^sq^{-b}c_{b,s}\gamma_{\sigma_u(j)-b}^{b-s+t-1}X^bH_{u,j}+c_{s-t,s}X^{s-t}K_{u,j}$
where $c_{s,s}=1,$ $c_{b,s}=0$ if $b<0,$ and
$c_{b,s}=\zeta_s\zeta_{s-1}\ldots \zeta_{b+1}$ with
$\zeta_s=(s)_{\!q}q^{-s}(q^{-(2j+u-1)-s}-1)$ if $0\ppq b\ppq s-1.$\el
\bpfl The proof is by induction on $t,$ and we  write out the case $t=1$ here:  
\begin{eqnarray*}
a_{\sigma_u(j)-s-1}X^sK_{u,j}&=& q^{-1}
Xa_{\sigma_u(j)-s}X^{s-1}K_{u,j}
-q^{-1}(1-q^{2j-2\rep{2j+u-1}-2s+u})X^{s-1}K_{u,j}\\
&=& q^{-1}
Xa_{\sigma_u(j)-s}X^{s-1}K_{u,j}
-(q^{-1}-q^{-2j-u+1-2s})X^{s-1}K_{u,j}\\
&=&q^{-1}\left(q^{-1}X^2a_{\sigma_u(j)-s+1}X^{s-2}K_{u,j}-(q^{-1}-q^{-2j-u+1-2s+2})X^{s-1}K_{u,j}\right)\\&&
-(q^{-1}-q^{-2j-u+1-2s})X^{s-1}K_{u,j}\\
&=&
q^{-2}X^2a_{\sigma_u(j)-s+1}X^{s-2}K_{u,j}\\&&-(q^{-1}+q^{-2}-q^{-2j-u+1-2s}-q^{-2j-u+1-2s+1})X^{s-1}K_{u,j}\\
&=& \ldots \\
&=& q^{-s}X^sa_jK_{u,j}+\zeta_s X^{s-1}K_{u,j}.
\end{eqnarray*}
\epfl

\bl\label{independence} For $0\ppq s \ppq d-\rep{2j+u-1}^--1,$ the elements $X^sK_{u,j}$ and
$X^s\tilde{F}_{u,j}$ are linearly independent (for other values of $s,$ we have $X^s\tilde{F}_{u,j}=0$).
\el

\bpfl Assume that $\alpha X^sK_{u,j}+\beta X^s \tilde{F}_{u,j}=0$ with
$0\ppq s \ppq d-\rep{2j+u-1}^--1.$ Multiply by
$\gamma_{\sigma_u(j)-1-s}^{s+1}$ using Lemma \ref{commut} and Lemma
\ref{R1}: then $\gamma_{\sigma_u(j)-1-s}^{s+1}X^s\tilde{F}_{u,j}$ is a
 multiple of
 $(q^{-2j+s+\rep{2j+u-1}-u+2-(s+1)}-1)X^{\rep{2j+u-1}}F_{u,j}$ which is zero. 

On the other hand,  $\gamma_{\sigma_u(j)-1-s}^{s+1}X^sK_{u,j}$ is
equal to
$\sum_{b=0}^sq^{-b}c_{b,s}\gamma_{\sigma_u(j)-b}^{b}X^bH_{u,j}.$ Now
$\zeta_p=0$ if and only if $p\equiv -2j-u+1;$ but if  $0\ppq b +1\ppq
p\ppq  s \ppq d-\rep{2j+u-1}^--1,$ we have  $\zeta_p\neq 0$ and
therefore $c_{b,s}\neq 0$ for all $b,s$ with $0\ppq b+1 \ppq s \ppq d-\rep{2j+u-1}^--1.$ 

So multiplying the identity  $\alpha X^sK_{u,j}+\beta X^s
\tilde{F}_{u,j}=0$ by $\gamma_{\sigma_u(j)-1-s}^{s+1}$ gives a nonzero
multiple of $\alpha \gamma_{\sigma_u(j)-1-s}^{s+1}X^sK_{u,j}$ with
$\gamma_{\sigma_u(j)-1-s}^{s+1}X^sK_{u,j}$ nonzero,  so $\alpha=0$ and
therefore $\beta=0.$\epfl

\bpfp {\bf \ref{projective}:} We  apply Lemma \ref{commut} with $t=d-1=s$ to see  that 
 $\gamma_{\sigma_u(j)-d}^{d-1}X^{d-1}K_{u,j}$ is a nonzero
 multiple of $\tilde{F}_{u,j}:$ if $b\ppq d-\rep{2j+u-1}^--1$ then
 $c_{b,d-1}=0;$ if $b \pgq d-\rep{2j+u-1}^-+1$ then
 $\gamma_{\sigma_u(j)-b}^{b-1}X^bH_{u,j}=0;$   finally, if
 $b=d-\rep{2j+u-1}^-, $ then
 $\gamma_{\sigma_u(j)-d}^{d-1}X^{d-1}K_{u,j}$ is a nonzero multiple of
  $X^{\rep{2j+u-1}}\gamma_j^{d-1}E_{u,j}.$
In particular, 
 $X^sK_{u,j}\neq 0$ for all $s=0,\ldots,d-1.$ 
The rest of the
 structure follows from this and Lemma \ref{independence}.
  \epfp

We can easily find all the  submodules of $\Gamma_uK_{u,j},$ and therefore:

\bdp When $\nproj,$ the module $\Gamma_uK_{u,j}$ has exactly two composition series:
$$\Gamma_uK_{u,j}\supset \Gamma_uF_{u,j}+\Gamma_u\tilde{D}_{u,j} \supset
\Gamma_u F_{u,j} \supset \Gamma_u \tilde{F}_{u,j}\supset 0$$ and $$\Gamma_uK_{u,j}\supset \Gamma_uF_{u,j}+\Gamma_u\tilde{D}_{u,j} \supset
\Gamma_u \tilde{D}_{u,j} \supset \Gamma_u \tilde{F}_{u,j}\supset 0.$$

Define $L(u,j):=\frac{\Gamma_uK_{u,j}}{\Gamma_uF_{u,j}+\Gamma_u\tilde{D}_{u,j}};$
 this is a simple module of dimension $d-\rep{2j+u-1}^-,$ and the
 composition factors of the composition series above are $L(u,j),$
 $L(u,\sigma_u(j)),$ $L(u,\sigma_u^{-1}(j)),$ and $L(u,j).$
 \edp 

\bpf The proof is straightforward using Proposition \ref{projective}
and the Proposition \ref{simples} which follows. \epf

\bp\label{simples} Let $S$ be a simple module. Then $S$ is isomorphic
to $L(u,j)$ if and only if the three following properties hold: 

\begin{enumerate}[(a)]
\item Dimensions: $\dim S=\dim L(u,j).$
\item Action of the central idempotents: $E_u$ acts as identity on
  $S$, and $E_v$ acts as zero on $S$ if $v\neq u$. 
\item Let $Y$ be the generator of $S$ which is in the kernel of the
  action of $X$ (this is well-defined up to a scalar and corresponds
  to $\tilde{H}_{v,j}$). Then the vertex $e_j$ acts as identity on
  $Y$, and the other vertices act as zero. 
\end{enumerate}
 
Note that the action of $G$ on this same element $Y$ is multiplication by $q^{j+u}.$

Moreover, $\dim L(u,\sigma_u^t(j))=\begin{cases} d-\rep{2j+u-1}^- & \mbox { if $t$ is even}\\\rep{2j+u-1} & \mbox { if $t$ is odd.} \end{cases}$ 

\ep

To summarize:

\bp If $\nproj,$ then $\Gamma_uK_{u,j}=\bigoplus
_{h=0}^{d-1}\left(kX^hF_{u,j}\oplus kX^hK_{u,j}\right)$ as a vector
space, it has dimension $2d, $ and it is an
indecomposable $\Gamma_u$-module  with simple top
and simple socle, both isomorphic to $L(u,j).$

If $\proj,$ then $L(u,j):=\Gamma_uK_{u,j}=\bigoplus
_{h=0}^{d-1}kX^hF_{u,j}$ as a vector space, it has dimension $d,$ and
it is a simple $\Gamma_u$-module.\ep

To decompose $\Gamma_u$ into a sum of indecomposable modules, we find modules isomorphic to the $\Gamma_uK_{u,j}$ inside $\Gamma_u:$

\bl If $0\ppq h\ppq d-\rep{2j+u-1}^--1,$ then right multiplication by
$X^h$ induces an isomorphism
$\Gamma_uK_{u,j}\stackrel{\sim}{\rightarrow}\Gamma_uK_{u,j}X^h$ of $\Gamma_u$-modules. \el

\bpfl  We must show that right multiplication by $X^h$ maps
$\Gamma_uK_{u,j}$ injectively in $\Gamma_uK_{u,j}X^h$, and it is
enough to do this for $h$ maximal. For this, we only need to check that
$\tilde{H}_{u,j}X^h\neq 0$ and $\tilde{D}_{u,j}X^h\neq 0.$ Since
$\tilde{H}_{u,j}X^h\neq 0$  implies $\tilde{D}_{u,j}X^h\neq 0,$ and since $h\ppq d-\rep{2j+u-1}^--1,$ we only need
to consider $\tilde{H}_{u,j}X^{d-\rep{2j+u-1}^--1}.$

To compute this, we need a relation similar to that in  Lemma \ref{R1}, which
is proved in the same way: we have
$$X^{d-1}\gamma_{j}^mE_{u,j}X^b=q^{-\frac{b(2m-b+1)}{2}}\frac{(m)!_{\!q}}{(m-b)!_{\!q}}\prod_{t=1}^b
(q^{2j+m+u+t}-1)\; X^{d-1}\gamma_{j+b}^{m-b}E_{u,j+b}$$ if $b<m$ and $X^{d-1}\gamma_{j}^mE_{u,j}X^b=0$ otherwise. Using this relation, we see that $\tilde{H}_{u,j}X^{d-\rep{2j+u-1}^--1}$
is a nonzero multiple of
$\prod_{t=1}^{d-\rep{2j+u-1}^--1}(q^{2j-1+u+t}-1)X^{d-1}\gamma_{j+d-\rep{2j+u-1}^--1}^{\rep{2j+u-1}^-}E_{u,-j-u}$
which is nonzero.
\epfl

We can now decompose $\Gamma_u$ entirely  into a sum of indecomposable
$\Gamma_u$-modules:

\bt $\Gamma_u$ decomposes into a direct sum of indecomposable modules
in the following way: $$\Gamma_u=\bigoplus_{j\in\zz_n}\bigoplus_{h=0}^{d-\rep{2j+u-1}^--1}\Gamma_uK_{u,j}X^h.$$ \et

\bpf We first prove that the sum is direct: the sums over $h$ (for $j$
fixed) are direct
because the  summands are in different right $G$-eigenspaces
($K_{u,j}X^hG=K_{u,j}X^hGE_{u,j+h}=q^{u+j+h}K_{u,j}X^h$).  The outer
sum is also direct because the summands have non-isomorphic socles:
the socle of $\bigoplus_{h=0}^{d-\rep{2j+u-1}^--1}\Gamma_uK_{u,j}X^h$ is
$\bigoplus_{h=0}^{d-\rep{2j+u-1}^--1}L(u,j)X^h$ with $L(u,j)X^h\cong
L(u,j).$

Equality follows from dimension counting.\epf

\bc\label{projectivesummary} Set $P(u,j)=\Gamma_uK_{u,j}$ for all $u,j.$ The modules $P(u,j)$
are projective, and they represent the different
isomorphism classes of projective $\dland$-modules when $u$ and $j$
vary in $\zz_n.$

When
$\nproj,$ their structure is
$$\xymatrix@C=.4cm@=.2cm{&L(u,j)\ar@{-}[dl]\ar@{-}[dr]\\
 L(u,\sigma_u^{-1}(j))&&L(u,\sigma_u(j))\\&L(u,j)\ar@{-}[ul]\ar@{-}[ur]}$$ 
and when $\proj,$ $P(u,j)=L(u,j)$ is simple of dimension $d.$

Moreover, the $L(u,j)$ represent all the isomorphism classes of simple
$\dland$-modules when $u$ and $j$ vary in $\zz_n$. Those of dimension $d$ are
also projective, and there are $\frac{n^2}{d}$ projective simples.
\ec

The simple modules were characterised in \cite{Rd2}, and when $n=d,$ the simple modules have been described in
 \cite{chen2}. 

We can now decompose each $\Gamma_u$, and therefore $\dland$, into
blocks: 

\bt\label{blocks} The algebras $\Gamma_u$ decompose into blocks as follows: if
$j_{u,1},\ldots,j_{u,r_u}$ are the representatives of the orbits of
$\sigma_u$ in $\zz_n,$ then  $\Gamma_u=\bigoplus_{{i=1}}^{r_u} \Bb_{{u,i}}$ where $$\Bb_{{u,i}}=\bigoplus_t \bigoplus_{h=0}^{d-\rep{2j+u-1}^--1}
P(u,\sigma^t(j_{u,i}))X^h,$$ where $t$ ranges from $0$ to
$\frac{2n}{d}-1$ if $2j_{u,i}+u-1 \nequiv 0$ and $t=0$ if $2j_{u,i}+u-1 \equiv 0.$

It then follows that the quiver of $\dland$ has $\frac{n^2}{d}$
isolated vertices which correspond to the simple projective modules, and
$\frac{n(d-1)}{2}$ copies of the quiver 
$$
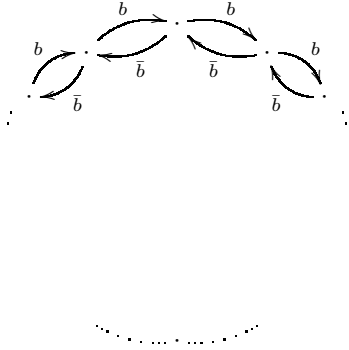
\xymatrix@=.01cm{
&&&&&&&&&&&\cdot\ar@/^.5pc/[rrrrrd]^{b}\ar@/^.5pc/[llllld]^{\bar{b}}\\
&&&&&&\cdot\ar@/^.5pc/[rrrrru]^{b}\ar@/^.5pc/[llldd]^{\bar{b}}&&&&&&&&&&\cdot\ar@/^.5pc/[rrrdd]^{b}\ar@/^.5pc/[lllllu]^{\bar{b}}\\\\
&&&\cdot\ar@/^.5pc/[rrruu]^{b}\ar@{.}@/_.3pc/[ldd]&&&&&&&&&&&&&&&&\cdot\ar@/^.5pc/[llluu]^{\bar{b}}\ar@{.}@/^.3pc/[rdd]\\\\
&&&&&&&&&&&&&&&&&&&&&&&\\
\\
\\
\\\\\\\\\\\\\\\\
&&&&&&&&&&&&&&&&\\
&&&&&&&&&&&\cdot\ar@/_.3pc/@{.}[rrrrru]\ar@/^.3pc/@{.}[lllllu]
}$$ with $\frac{2n}{d}$ vertices and $\frac{4n}{d}$ arrows. The
relations on this quiver are $bb$, $\bar{b}\bar{b}$ and
$b\bar{b}-\bar{b}b$ (there are $\frac{6n}{d}$ relations on each of
these quivers). The vertices in this quiver correspond to the simple modules $L(u,j),$
$L(u,\sigma_u(j)),$ $L(u,\sigma_u^2(j)),$ $\ldots,$ $L(u,\sigma_u^{\frac{2n}{d}-1}(j)).$\et

\bpf For the general principle of presenting a basic algebra by quiver
and relations, see \cite[Section II.5 and III.1 Theorem 1.9]{ARS}. 
The arrows of the quiver correspond to generators of $\rad(\dland)/\rad^2(\dland)$.
From the structure of the
indecomposable projective modules we
know that for each $p$ in $\zz_{\frac{2n}{d}}$, there is one
arrow which we call $b_p$ from
 the vertex $\epsilon_p$ -- corresponding to the simple module
 $L(u,\sigma_u^p(i))$ -- to the vertex  $\epsilon_{p+1}$, and one arrow,
called $\overline{b}_p$, from $\epsilon_{p+1}$ to $\epsilon_p$.
The zero relations follow easily, and moreover, for each $p$,
there is a non-zero scalar
$c_p$ with $c_p\overline{b}_pb_p = b_{p-1}\overline{b}_{p-1}$.

\medskip

Starting at $p=1$, we replace $c_pc_{p-1}\ldots c_1\overline{b}_p$ by $\overline{b}_p$, for $p=1,
2,\cdots, \frac{2n}{d}-1$, so that the relations become $$
\begin{array}{ll}
\overline{b}_pb_p = b_{p-1}\overline{b}_{p-1} & \mbox{ for $1\ppq p
  \ppq \frac{2n}{d}-1$ }\\
c_0\ldots c_{\frac{2n}{d}-1}\overline{b}_0b_0=b_{-1}\overline{b}_{-1}.

\end{array}
$$

Since the algebra is symmetric (see the Introduction), the scalar in the final relation is 1:
to see this, take a symmetrising form $\psi:\mathbb{B}_{u,i}\to k$; then
$$\psi(c_0\ldots c_{\frac{2n}{d}-1}\overline{b}_0b_0) = \psi(b_{-1}\overline{b}_{-1}) 
=\psi(\overline{b}_{-1}b_{-1})
=\psi(b_0\overline{b}_0) = \psi(\overline{b}_0b_0).
$$
Hence $(c_0\ldots c_{\frac{2n}{d}-1}-1)\overline{b}_0b_0$ lies in the
kernel of $\psi$, which spans a 1-dimensional (left) ideal
of the algebra, and since $\psi$ is non-singular, this must be zero. This shows that $c_0\ldots c_{\frac{2n}{d}-1}=1$.

\epf


\br\label{specialbiserialandrepresentationtype} It is now easy to see from the quiver and the structure of the projectives that the algebra $\dland$ is
special biserial (see \cite[II.1]{erdmann}). 

Note also that it follows from \cite[II.3.1]{erdmann} that $\dland$ is therefore tame or of finite type (see \cite[p111]{ARS} and \cite[4.4]{benson} for the definitions). We shall see in Section \ref{sectionclassification} and in the Appendix that it is in fact tame.\er

\br  It is known that the finite-dimensional quotients of
$U_q(\frak{sl}_2)$ studied by R. Suter \cite{suter}, J. Xiao
\cite{xiao} and M. Patra \cite{patra} are quotients of some of these Drinfel'd doubles (for
particular choices of $d$). Therefore their results can be recovered
from the study of the Drinfel'd doubles $\dland$.

\er

\br 
If we fix a block $\mathbb{B}$ of $\dland$ which is not simple, we can
prove that it is a Koszul algebra, using the results in \cite[Section 3]{GM} (by looking at  minimal projective
resolutions of the simple modules for  $\mathbb{B}$). Therefore, by
\cite[Theorem 6.1 and Section 10]{GM}, its Koszul dual
$\Ext^*_\mathbb{B}(\mathbb{B}/\rad(\mathbb{B}),\mathbb{B}/\rad(\mathbb{B}))$
is given by $k\mathcal{Q}^{op}/I^\perp$, where $\mathcal{Q}$ is the quiver of
$\mathbb{B}$ and $I^\perp$ is generated by the relations
$b\overline{b}+\overline{b}b$ for all the arrows $b,\overline{b}$ in
the quiver.

In fact, we see that
$\Ext^*_\mathbb{B}(\mathbb{B}/\rad(\mathbb{B}),\mathbb{B}/\rad(\mathbb{B}))$
is the preprojective algebra associated to an (unoriented) cycle. In
\cite[Theorem 7.2]{G}, it was shown that such an algebra is Koszul, giving another
proof that $\mathbb{B}$ is a Koszul algebra.
\er


\section{Classification of the representations of $\dland$}\label{sectionclassification}

We can now classify all the indecomposable representations of $\dland$ and describe
its Auslander-Reiten quiver. For definitions of {\bf Auslander-Reiten sequences} (also called {\bf almost split sequences}) and {\bf Auslander-Reiten quivers}, as well as descriptions of some components of Auslander-Reiten quivers, see \cite[V.1 and VII.1]{ARS} and \cite[I.7 and I.8]{erdmann}. In this section we describe the
representations and the Auslander-Reiten quiver without
proofs; these are outlined in the Appendix. Note that this description shows that the algebra $\dland$ is tame (see Remark \ref{specialbiserialandrepresentationtype}). We then determine the quantum
dimension of these representations and study related properties which
we use in Section \ref{tensor} when we calculate the tensor products of representations.

\subsection{Description of the indecomposable modules}\label{desc}

\begin{enumerate}[{\bf (I)}]
\item\label{descodd}{\bf String modules of odd length:} the indecomposable modules of
  odd length are syzygies of simple modules, that is, of the form $\Omega^k(L(u,i))$ for some $k$ in $\zz$
  and for some simple module $L(u,i).$
  We have $$\abs{{\rm length}(\Top (\Omega^k(L(u,i))))-{\rm
      length}(\soc (\Omega^k(L(u,i))))}=1.$$ These modules are not periodic.

\item\label{descstring} {\bf String modules of even length:} fix a block $\Bb_{u,i}.$
  For each $0\ppq p\ppq \frac{2n}{d}-1$ and for each $\ell \pgq 1,$ there are two indecomposable
  modules of length $2\ell$ which we call $M^{\pm }_{2\ell}(u,\sigma_u^p(i)):$ 
\begin{enumerate}[$\bullet$]
\item The module $M^{+}_{2\ell}(u,\sigma_u^p(i))$ has top composition factors
  $L(u,\sigma_u^p(i)),$
  $L(u,\sigma_u^{p+2}(i)),$ $\ldots,$ $L(u,\sigma_u^{p+2(\ell-1)}(i))$ and
    socle composition factors $L(u,\sigma_u^{p+1}(i)),$ $
  L(u,\sigma_u^{p+3}(i)),$ $\ldots,$\\
  $L(u,\sigma_u^{p+2(\ell-1)+1}(i))$:
$$\def\objectstyle{\scriptstyle}\xymatrix@=.05cm{L(u,\sigma_u^p(i)) \ar@{-}[dr]  &  &
    L(u,\sigma_u^{p+2}(i))\ar@{-}[dr]\ar@{-}[dl] &  & 
  L(u,\sigma_u^{p+2(\ell-1)}(i))\ar@{-}[dr]\ar@{-}[dl]&&              
 \\&L(u,\sigma_u^{p+1}(i)) && \cdots &&
L(u,\sigma_u^{p+2(\ell-1)+1}(i)) &&     
}$$ The lines joining the simple modules are given by multiplication
by the appropriate $b$-arrow or $\bar{b}$-arrow (in the case $n=d$,
when there is an ambiguity, the first line is multiplication by
$\gamma^{\dim L(u,\sigma_u^p(i))}$, the next one is multiplication by
a scalar multiple of $X^{d-\dim L(u,\sigma_u^p(i))}$, and so on, up to scalars).
\item The module $M^{-}_{2\ell}(u,\sigma_u^p(i))$ has top composition factors
  $L(u,\sigma_u^p(i)),$ $
  L(u,\sigma_u^{p-2}(i)),$ $\ldots,$ $L(u,\sigma_u^{p-2(\ell-1)}(i))$ and
    socle composition factors $L(u,\sigma_u^{p-1}(i)),$ $
  L(u,\sigma_u^{p-3}(i)),$ $\ldots,$\\
  $L(u,\sigma_u^{p-2(\ell-1)+1}(i))$:
$$\def\objectstyle{\scriptstyle}\xymatrix@=.05cm{ & L(u,\sigma_u^{p-2(\ell-1)}(i))\ar@{-}[dr]\ar@{-}[dl]  &&  L(u,\sigma_u^{p-2}(i)) \ar@{-}[dr] \ar@{-}[dl] &&
L(u,\sigma_u^p(i)) \ar@{-}[dl]  
\\L(u,\sigma_u^{p-2(\ell-1)+1}(i)) &&   \ldots && L(u,\sigma_u^{p-1}(i))
} $$ As for the other string modules, the lines represent
multiplication by an appropriate $b$ or $\bar{b}$ arrow, and in the
case $n=d$ the first one from the left is multiplication by a power of $X$ and so on.
\end{enumerate} In both cases, indices are taken modulo $\frac{2n}{d}.$

These modules are periodic of period $\frac{2n}{d}.$

\item\label{descband} {\bf Band modules (even length):} fix a block $\Bb_{u,i}$. For
  each $\lambda\neq 0$ in $k,$ and for each $\ell\pgq 1$, there are two
  indecomposable modules of length $\frac{2n}{d}\ell$, which we denote
  by $C^{\ell\pm}_\lambda(u,i)$. They are
  defined as follows:

\begin{enumerate}[$\bullet$]

\item Let $V$ be an $\ell$-dimensional vector space. Then
  $C^{\ell+}_\lambda(u,i)$ has underlying space
  $C^+=\bigoplus_{p=0}^{\frac{2n}{d}-1}C^+_p$ with $C^+_p=V$ for all $p.$ The
  action of the idempotents $\epsilon_p$ is such that
  $\epsilon_pC^+=C^+_p.$ The action of the arrows $\bar{b}_{2p}$ and
  $b_{2p+1}$ is zero. The action of the arrows $\bar{b}_{2p+1}$ is the
  identity of $V.$ The action of the arrows $b_{2p}$ with $p\neq 0$ is
  also the identity. Finally, the action of $b_0$ is given by the
  indecomposable Jordan matrix $J_\ell(\lambda).$

Note that $\soc(C^+)=\rad(C^+)=\bigoplus_{p}{\epsilon_{2p+1}C^+}$ and that
$C^+/\rad(C^+)=\bigoplus_{p}{\epsilon_{2p}C^+}$.

\item The module $C^{\ell-}_\lambda(u,i)$ is defined similarly,
  interchanging $b$'s and $\bar{b}$'s. Note that if $C^-=\bigoplus_{p=0}^{\frac{2n}{d}-1}C_p^-$ is the
  underlying vector space of $C^{\ell-}_\lambda(u,i)$, the Jordan matrix
  $J_\ell(\lambda)$ occurs as $\bar{b}_0$ from $C^-_1$ to
  $C^-_{0},$ and that we have
  $\soc(C^-)=\bigoplus_{p}{\epsilon_{2p}C^-}$ and
  $C^-/\rad(C^-)=\bigoplus_{p}{\epsilon_{2p+1}C^-}.$
 
\end{enumerate}

These modules are periodic of period 2.

\end{enumerate}


\subsection{Description of the Auslander-Reiten components}\label{arcomponents}

\begin{enumerate}[{\bf (I)}]
\item\label{arodd} {\bf Components with indecomposable modules of odd length:} there are
  two components in the Auslander-Reiten quiver for each block $\Bb_{u,i}$, one of which contains the
  simple modules $L(u,\sigma_u^{2p}(i))$ and the other containing the
  simple modules $L(u,\sigma_u^{2p+1}(i))$ for all $0\ppq p\ppq
  \frac{2n}{d}-1.$ They are of tree class
  $\tilde{A}_{\frac{2n}{d}-1},$ infinite in all directions, but
  we identify along horizontal lines so that the component lies on an infinite cylinder:

$$\def\objectstyle{\scriptstyle}
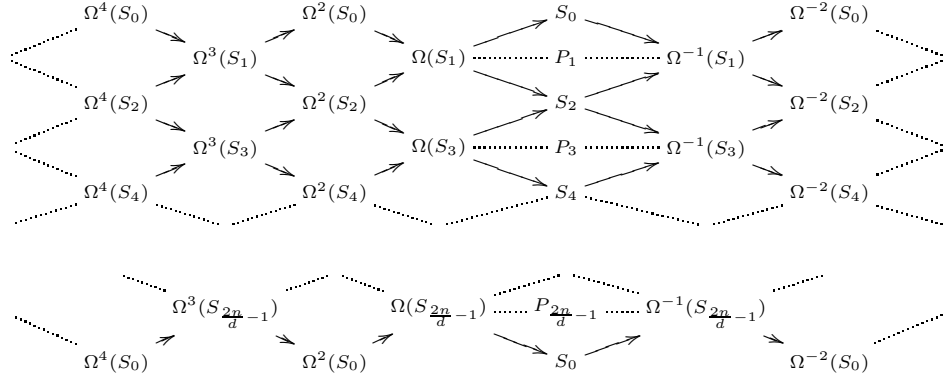
\xymatrix@=.1cm{&\hspace{1.5cm} &\Omega^4(S_0)\ar[dr]\ar@{.}[dl]&& \Omega^2(S_0)\ar[dr]& & &S_0\ar[drr] &&& \Omega^{-2}(S_0)\ar@{.}[dr] &  \hspace{1.5cm} &\\&
&\ &\Omega^3(S_1)\ar[dr]\ar[ur] &&\Omega(S_1)\ar[drr]\ar[urr]\ar@{.}[rr] &&P_1\ar@{.}[rr] &&\Omega^{-1}(S_1)\ar[dr]\ar[ur]& & &\\&
\ &  \Omega^4(S_2)\ar@{.}[dl]\ar@{.}[ul]\ar[dr]\ar[ur] &&
\Omega^2(S_2)\ar[dr]\ar[ur] &&&S_2\ar[drr]\ar[urr] &&&\Omega^{-2}(S_2)\ar@{.}[dr]\ar@{.}[ur]& & &\\&
\ &&\Omega^3(S_3)\ar[dr]\ar[ur] &&\Omega(S_3)\ar[drr]\ar[urr]\ar@{.}[rr] &&P_3\ar@{.}[rr] &&\Omega^{-1}(S_3)\ar[dr]\ar[ur]& & &\\&
\ &\Omega^4(S_4)\ar@{.}[dl]\ar@{.}[ul]\ar[ur]\ar@{.}[dr]&& \Omega^2(S_4)\ar[ur]\ar@{.}[dl]\ar@{.}[dr]  & &&
S_4\ar[urr]\ar@{.}[drr] \ar@{.}[dll]&&& \Omega^{-2}(S_4)\ar@{.}[dl]\ar@{.}[dr]\ar@{.}[ur] & &\\&
\ &&&
&&&&&&&& &\\&
&\ &&&&&&&&&&& &\\&
\ &&&&&&&& &&& &\\&
\ &&\Omega^3(S_{{\frac{2n}{d}}-1})\ar@{.}[ul]\ar@{.}[ur]\ar[dr] &&\Omega(S_{{\frac{2n}{d}}-1})\ar[drr]\ar@{.}[rr]\ar@{.}[ul]\ar@{.}[urr] &&P_{{\frac{2n}{d}}-1}\ar@{.}[rr] &&\Omega^{-1}(S_{{\frac{2n}{d}}-1})\ar[dr]\ar@{.}[ull]\ar@{.}[ur]& & &\\&
\ &\Omega^4(S_0)\ar@{.}[ul]\ar[ur]&& \Omega^2(S_0)\ar[ur]  &&& 
S_0\ar[urr] &&& \Omega^{-2}(S_0)\ar@{.}[ur] & &
}$$ where $S_p$ represents $L(u,\sigma_u^p(i))$, and $P_p$ is the
projective $P(u,\sigma_u^p(i)).$

\item\label{arstring} {\bf Components with string modules of even length:} there are
  four components in the Auslander-Reiten quiver  for each block $\Bb_{u,i},$ one  which contains
the modules $M^{+}_{2\ell}(u,\sigma_u^{2p}(i))$ for all $p$ and $\ell,$ one which
contains all the modules $M^{+}_{2\ell}(u,\sigma_u^{2p+1}(i))$ for all $p$ and
$\ell,$ one which contains the modules $M^{-}_{2\ell}(u,\sigma_u^{2p}(i))$ for all $p$ and $\ell,$
and the final one contains all the modules $M^{-}_{2\ell}(u,\sigma_u^{2p+1}(i))$ for all $p$ and
$\ell.$ They are tubes of rank $\frac{n}{d},$ and the modules of
length $2\ell$ form the $\ell^{ th}$ row.

The Auslander-Reiten sequences are: $$
0\rightarrow M^+_{2\ell}(u,i-d)\rightarrow 
\begin{array}{c}
M^+_{2\ell+2}(u,i-d)\\\oplus\\
M^+_{2\ell-2}(u,i)
\end{array}\rightarrow M^+_{2\ell}(u,i)\rightarrow 0
$$ and 
$$
0\rightarrow M^-_{2\ell}(u,i+d)\rightarrow 
\begin{array}{c}
M^-_{2\ell+2}(u,i+d)\\\oplus\\
M^-_{2\ell-2}(u,i)
\end{array}\rightarrow M^-_{2\ell}(u,i)\rightarrow 0
$$ (where $M^\pm_0(u,i):=0$).

\item\label{arband}  {\bf Components with band modules:} fix a block
  $\Bb_{u,i}$. There are two components in the Auslander-Reiten quiver  for
each nonzero parameter $\lambda$ in $k$, one which contains  the
modules $C^{\ell +}_\lambda(u,i)$ for all $\ell$ and the other which contains the
modules $C^{\ell -}_\lambda(u,i)$  for all $\ell$. They are tubes of
rank one.\end{enumerate}

The Auslander-Reiten sequences are: $$ 0\rightarrow C_\lambda^{\ell +}(u,i)\rightarrow 
\begin{array}{c}
C_\lambda^{\ell+1, +}(u,i)\\\oplus\\C_\lambda^{\ell-1,
  +}(u,i)
\end{array}\rightarrow C_\lambda^{\ell +}(u,i)\rightarrow 0
$$ and 
$$ 0\rightarrow C_\lambda^{\ell -}(u,i)\rightarrow 
\begin{array}{c}
C_\lambda^{\ell+1, -}(u,i)\\\oplus\\C_\lambda^{\ell-1,
  -}(u,i)
\end{array}\rightarrow C_\lambda^{\ell -}(u,i)\rightarrow 0
$$ where $C_\lambda^{0\pm}(u,i):=0.$

\subsection{Splitting trace modules}\label{sectionsplittingtracemodules}

In this section, we determine the quantum dimensions of the modules
described above. This concept is defined in \cite{charipressley} for
ribbon algebras (note that the Drinfel'd double of a Hopf algebra
naturally gives rise to a ribbon algebra and that $S(x)=GxG^{-1}$ for
all $x\in \dland$):

\bd\cite[Section 4.2.C]{charipressley} Let $\rho:\dland\rightarrow \End_k(V)$ be a
representation of $\dland.$ Let $f:V\rightarrow V$  be a linear map
and let $\tr$ be the usual trace.
\begin{enumerate}[$\bullet$]
\item The quantum trace of $f$ is $\qtr(f)=\tr(\rho(G)f).$
\item The quantum dimension of $V$ is $\qdim(V)=\qtr(\id_V).$
\end{enumerate}

\ed

In Section \ref{tensor}, we shall consider tensor products of modules
over $k,$ and in this context, {\bf splitting trace modules,} that
is, modules $M$ such that the trivial module $L(0,0)$ is a direct
summand in $\End_{\dland}(M),$ are useful in view of the following
proposition:

\bp\label{arandstm} Let $M$ be an indecomposable
module, and let  $\A(M)$ be its
Auslander-Reiten sequence. 

\begin{enumerate}[(i)]
\item\label{arstm} 
Suppose $N$ is a splitting trace module, and 
  $M\ot N= \bigoplus _i A_i\oplus \bigoplus_j B_j\oplus P$, with the $A_i$
  non-projective indecomposable  splitting trace modules, the $B_j$
  non-projective  indecomposable and not splitting trace modules, and $P$
  projective. Then $\A(M)\ot N$ is the direct sum of
  $\bigoplus_i \A(A_i)$ and of a split exact sequence (equal to
  $0\rightarrow\bigoplus _j \Omega^2(B_j)\oplus Q \rightarrow \bigoplus
  _j \Omega^2(B_j)\oplus \bigoplus_j B_j \oplus Q \oplus P \rightarrow
  \bigoplus_j B_j  \oplus P \rightarrow 0 $ for some projective $Q$).
\item\label{arnotstm} If $N$ is not a splitting trace module,
then $\A(M)\ot N$ is split exact.
\end{enumerate}
\ep

This proposition is proved using the same methods as \cite[2.5.9]{H},
\cite[2.6]{AC}, \cite[3.1 and 3.8]{GMS} and \cite[1.5]{Ka}, replacing
the natural isomorphism $M\cong M^{**}$ by the isomorphism of
$\dland$-modules given by $m\mapsto \rep{-,Gm}.$

Splitting trace modules are related to the quantum trace as follows:

\bp\cite[2.5.9]{H} A module $M$ is a splitting trace module if and
only if there exists an endomorphism
$f\in\End_{\dland}(M)$ such that $\qtr(f)\neq 0$.
\ep 

\br Note that if $\qdim(M)\neq 0,$ then $M$ is a splitting trace module. \er

\bp\label{stm} Let $M=\Omega^\ell(L(u,i))$ be a module of odd length. Then
$\displaystyle{\qdim(M)=(-1)^\ell q^{i+u}\frac{1-q^{\dim(L(u,i))}}{1-q}}$. Therefore if  $M$ is not projective,
then it  is a splitting trace module.

\ep

\bpf We prove it for $\ell$ non-negative and even; the other cases are similar: $$M=\ \def\objectstyle{\scriptstyle}\xymatrix@=.05cm{L(u,\sigma_u^{-\ell}(i)) \ar@{-}[dr]  &  &
    L(u,\sigma_u^{\ell-2}(i))\ar@{-}[dr]\ar@{-}[dl] &  & 
  L(u,i)\ar@{-}[dl]\ar@{-}[dr]&& L(u,\sigma_u^{\ell}(i))  \ar@{-}[dl]         
 \\&L(u,\sigma_u^{-\ell+1}(i))\ar@{-}[ru] && \cdots &&\cdots
   }.$$ Set $N=\dim L(u,i).$ 
For each composition factor, choose a basis obtained by taking a generator of this
simple which is in the kernel of the action of $X$, and then applying
arrows to this element. These bases can be chosen so that they are compatible with the actions of the $b$- and $\bar{b}$-arrows. This gives a basis for $M$, with respect to
which the action of $G$ can be described by the diagonal matrix with entries 
$
q^{i+u},$ $
 q^{i+u+1},$ $
 \ldots,$ $
 q^{i+u+N-1},$ $
 q^{\sigma_u(i)+u},$ $
 \ldots,$ $
 q^{\sigma_u(i)+u+d-N-1},$ $
 q^{i+u},$ $
 q^{i+u+1},$ $
 \ldots,$ $
 q^{i+u+N-1},$ $
 \ldots,$ $
 q^{i+u+N-1}$,
 so that $\qdim(M),$ which is the trace of
this matrix, is $$\ell
\left(\sum_{t=0}^{N-1}q^{i+u+t}+\sum_{s=0}^{d-N-1}q^{\sigma_u(i)+u+s}\right)+\sum_{t=0}^{N-1}q^{i+u+t}.$$
Since $\sigma_u(i)=i+N,$ setting $s=t-N$ gives:
$$\qdim(M)=\ell\sum_{t=0}^{d-1}q^{i+u+t}+\sum_{t=0}^{N-1}q^{i+u+t}=\ell
q^{i+u}\frac{1-q^d}{1-q}+q^{i+u}\frac{1-q^N}{1-q}=q^{i+u}\frac{1-q^N}{1-q}.$$ Therefore $M$ is a splitting trace module if, and only if, $N\neq d,$ that is, $M$ is not projective.

\epf

\bp\label{notstm} Let $M$ be an indecomposable module of even
length. Then $M$ is not a splitting trace module and $\qdim(M)=0$.
\ep

\bpf
Take
$f$ in $\End_{\dland}(M).$ We choose a basis for $M$ as in the proof of Proposition \ref{stm}. Since $f$ is a
homomorphism of $\dland$-modules, the matrix of $f$ with respect to
this basis is diagonal. By considering the actions of the arrows $b_p$
and $\bar{b}_p$ which connect the composition factors of $M$, we can
see that $f$ is of the form $\mu \id$ for some $\mu\in k$. Therefore
$\qtr(f)=\mu\qtr(\id)=\mu\qdim(M).$ We can then calculate this as in the proof
of Proposition \ref{stm}, to see that $\qdim(M)=0$. Therefore
$\qtr(f)=0$ for all $f\in\End_{\dland}(M)$ and so $M$ is
not a splitting trace module.
\epf


\section{Tensor products of 
        $\dland$-modules}\label{tensor}

In this section, we determine the tensor products over the base field
$k$ of modules of types {\bf (I)} and {\bf (II)}. Recall that if $M$
and $N$ are modules over $\dland,$ then $M\ot N \cong N\ot M$ (where
$\ot=\ot_k$).

\subsection{Tensor product of simple $\dland$-modules}

The aim of this section is to prove the following theorem:

\bt\label{tensorsimples} Set $\varpi=\dim L(u,i)+\dim L(v,j)-(d+1),$ and define
$$\varsigma=\begin{cases}
  \frac{\varpi}{2} &\mbox{ if $\varpi$ is even }\\
  \frac{\varpi+1}{2} &\mbox{ if $\varpi$ is odd.} \end{cases}$$

Then we have the following decompositions:

\begin{enumerate}[(a)]
\item If $\varpi\ppq 0,$ $$L(u,i)\ot L(v,j)\cong
  \bigoplus_{\theta=0}^{\mathrm{min}\set{\dim L(u,i)-1;\dim
      L(v,j)-1}}L(u+v,i+j+\theta).$$ 
\item If $\varpi\pgq 0,$ $$L(u,i)\ot L(v,j)\cong  \bigoplus_{\theta=\varpi+1}^{\mathrm{min}\set{\dim L(u,i)-1;\dim
    L(v,j)-1}}L(u+v,i+j+\theta)\oplus \bigoplus_{\theta=\varsigma}^{\varpi}P(u+v,i+j+\theta).$$
\end{enumerate}  \et

\br  In the case $n=d,$ H-X. Chen gives this decomposition when
$\varpi\ppq 0$, and gives the socle of this tensor product when
$\varpi \pgq 0$
in \cite{chen2}.\er


\subsubsection{Preliminaries}

Recall (Proposition \ref{simples}) that a simple module is determined
up to isomorphism by its dimension, the action of the
idempotents $E_u$ and the action of the idempotents $e_i.$ This will
be important in order to determine the socle of the tensor product
$L(u,i)\ot L(v,j)$.

\bn \label{dimassumptions} To simplify notation, set $$\begin{cases}\dim
    L(u,i)=\alpha,\\  \dim L(v,j)=\beta,\\  \mbox{assume that
$\alpha\ppq \beta $ throughout.}\end{cases}$$\en

\bl The vector space $L(u,i)\ot L(v,j)$ has basis
$\set{z_{s,t}\,\mid\, 0\ppq s \ppq \alpha-1,\; 0 \ppq t \ppq \beta-1},$ where
$z_{s,t}=\gamma_i^s\vv{u,i}\ot \gamma_j^t\vv{v,j}.$ 
 \el

\bd We define a {\bf grading} on the vector space  $L(u,i)\ot L(v,j)$ by
setting $\gr(z_{s,t})=s+t.$\ed

\bl The socle of  $L(u,i)\ot L(v,j)$ is isomorphic to a direct sum of simple modules
of the form $L(u+v, -).$\el

\bpfl It is easy to see that $E_{u+v}z_{s,t}=z_{s,t}$ for all $s,t.$ \epfl
 
In order to find all the simple modules in the socle, we need to find
the kernel of the action of $X$: the generators
$\vv{w,\ell}$ of the simple modules are in the kernel of the action of $X$. Note that $X\cdot
\elt{w}{\ell}{s}=(s)_{\!q}(1-q^{2\ell+w-1+t})\elt{w}{\ell}{s-1}$, which is nonzero
if $s>0,$ and
zero if $s=0.$


\subsubsection{Kernel of $X$}

\bl The kernel of the action of $X$ on $\dland$
has dimension $n^2d.$
\el

\bpfl The Hopf algebra $\land^{*cop}$ is a Hopf subalgebra of $\dland,$ so by
a theorem of Nichols-Zoeller \cite[Theorem 7]{nichols-zoeller}, $\dland$ is
free as a module over $\land^{*cop},$ with rank $nd.$ 

Now $\land^{*cop}=\rep{G,X \,\mid \, G^n=1, X^d=0, GX=q^{-1}XG}$ so
the kernel of the action of $X$ on $\land^{*cop}$ has dimension $n$
with basis $\set{G^iX^{d-1}\,\mid\, i=0,\ldots,n-1}.$

Therefore the kernel of the action of $X$ on $\dland$ has dimension
$n\cdot nd.$  \epfl

\bl The kernel of $X$ is spanned by $\set{\vv{u,j}X^h,\tilde{D}_{u,j}X^h\,
  \mid \, u,j \in \zz_n, 0\ppq h\ppq d-\rep{2j+u-1}^-}$ (this set
contains repetitions). \el

\bpfl Note that all the elements $\vv{u,j}X^h$ and $\tilde{D}_{u,j}X^h$ are
in the kernel of the action of $X.$ 

We consider the case where $d$ is odd; the case $d$ even is similar,
but the simple projective modules are distributed differently.

  There are $n$ simple projective modules of
dimension $d$ in $\Gamma_u$ (there are $\frac{n}{d}$ isomorphism classes of
simple projective modules, and the number of isomorphic copies of each
one is $d$) and each one gives one element of the form
$\vv{u,j}X^h=\tilde{D}_{u,j}X^h.$ The other projective modules in $\Gamma_u$
have dimension $2d,$ so there are $\frac{\dim \Gamma_u - n\cdot
  d}{2d}=\frac{n(d-1)}{2}$ of them, and each one gives two elements in the
kernel of $X.$ Therefore we have $n\cdot (n+2\frac{n(d-1)}{2})=n^2d$ distinct
elements in  $\set{\vv{u,j}X^h,\tilde{D}_{u,j}X^h\,
  \mid \, u,j \in \zz_n, 0\ppq h\ppq d-\rep{2j+u-1}^-}$  which are in the kernel of $X$ and linearly
independent. We can visualise these elements in the basis of the action
of $X$ as in the picture below: $$
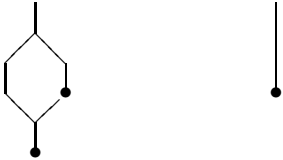
\xymatrix@=.4cm@M=0cm@!0{&\ar@{-}[d]&&&&&&&&  \ar@{-}[ddd]    \\&\ar@{-}[dl]\ar@{-}[dr]&&&&&&&&    \\\ar@{-}[d]&&\ar@{-}[d]&&&&&&&    \\\ar@{-}[dr]&&\bullet
  \ar@{-}[dl]&&&&&&&    \bullet  \\ &\ar@{-}[d] &&& \\&\bullet&&&&&&&&}$$

Note that there are $n^2d+n^2$ elements in $\set{\vv{u,j}X^h,\tilde{D}_{u,j}X^h\,
  \mid \, u,j \in \zz_n, 0\ppq h\ppq d-\rep{2j+u-1}^-}$, which is
consistent with the fact that the elements in the kernel of the action
of $X$ coming from simple projective modules appear
twice. \epfl

From now on, we consider  modules up to isomorphism,
so assume that $h=0.$

\bp\label{kernel}  The kernel of the action of $X$ on $L(u,i)\ot L(v,j)$ has
dimension $\alpha=\dim L(u,i),$ and is spanned by elements $x_\theta$ with
$\gr{x_\theta}=\theta$ and $0\ppq \theta \ppq \alpha-1$ (there is one
element (up to scalars) in the kernel for each degree).\ep

\bpf Consider the vector space $U_\theta$ spanned by the elements
$z_{s,t}$ with $\gr(z_{s,t})=\theta.$ The kernel of the action of $X$
on  $L(u,i)\ot L(v,j)$ is the span of the kernels of the
$\varphi_\theta$ where $\varphi_\theta:U_\theta\rightarrow
U_{\theta-1}$ is right multiplication by $X.$ It is easy to see that
$\varphi_\theta(z_{s,t})=\mu_{s,t}z_{s-1,t}+\nu_{s,t}z_{s,t-1}$ with
$\mu_{s,t}=0$ if and only if $s=0$ and $\nu_{s,t}=0$ if and only if
$t=0.$

First assume that $0\ppq \theta \ppq \alpha-1:$ then $\dim
U_\theta=\theta+1$ and $\dim U_{\theta-1}=\theta.$ The matrix of
$\varphi_\theta$ is of the form  $$
\left(\begin{array}{c|cccc}
\mu & \nu \\
&\mu&\ddots &0\\
&0&\ddots &\nu \\
&&&\mu&\nu
\end{array}\right)
$$ so $\varphi_\theta$ is onto; therefore its kernel has dimension 1 and is
spanned by an element $x_\theta$ with $\gr(x_\theta)=\theta.$

If $\theta\pgq \alpha,$ then $\dim U_{\theta-1}\pgq \dim U_{\theta}$
and the matrix of $\varphi_\theta$ is either square (if $\alpha\ppq
\theta\ppq \beta-1$) or rectangular with one more line than it has
columns. It contains a maximal square matrix 
of the form $$
\left(\begin{array}{cccc}
 \nu \\
\mu&\ddots &0\\
0&\ddots &\nu \\
&&\mu&\nu
\end{array}\right)
$$
so the rank of $\varphi_\theta$ is  the dimension of
$U_\theta,$ and $\varphi_\theta$ is therefore injective.  \epf

\br An element $x_\theta$ is either of the form $\vv{u+v,p}$, in
which case it generates a simple module, or of the form
$\tilde{D}_{u+v,p},$ in which case it generates a non-simple 
module which cannot be a summand in the socle of $L(u,i)\ot L(v,j)$.
\er

\br Since $e_{i+j+\theta}z_{s,\theta-s}=z_{s,\theta-s}$ for all $s,$ we also have
$e_{i+j+\theta}x_\theta=x_\theta.$ In particular, if $x_\theta$
generates a simple module, then this simple module must be $L(u+v,i+j+\theta).$ \er


\subsubsection{Proof of the decomposition when $\alpha+\beta\ppq d+1,$ and
  some general results}

\bp  Suppose $\varpi<0$ and $0\ppq \theta<\alpha$ or $\varpi\pgq 0$
and  $0\ppq\varpi< \theta<\alpha .$ Then $x_\theta$ generates a  simple
module of dimension $<d.$  The element $x_\theta$ for
$2\theta=\varpi$ generates a simple projective module if and only if $\varpi> 0$ is even or 
$\varpi=0$.

Therefore $L(u,i)\ot
L(v,j)$ contains $\bigoplus_{\theta=0}^{\alpha-1}L(u+v,i+j+\theta)$ if $\varpi<0$ and contains
$\bigoplus_{\theta=\varpi+1}^{\alpha-1}L(u+v,i+j+\theta)\oplus
L(u+v,i+j+\frac{\varpi}{2})$ if $\varpi\pgq 0$, the last term only occurring if
$\varpi$ is   even.  \ep

\bpf  If $x_\theta$ is in the kernel of $\gamma^{d-1}_{i+j+\theta},$
then it must generate a simple module of dimension $<d$ ({see}
Proposition \ref{projective}). So suppose for a contradiction that
$\gamma^{d-1}_{i+j+\theta}x_\theta\neq 0.$ Then
$\gamma^{d-1}_{i+j+\theta}z_{s,\theta-s}\neq 0$ for some $s.$ But
$\gamma^{d-1}_{i+j+\theta}z_{s,\theta-s}   =\sum_{r=0}^{d-1}\begin{pmatrix}d-1\\r\end{pmatrix}_{\!q}q^{r(\theta-s)}\elt{u}{i}{r+s}\ot
\elt{v}{j}{\theta+d-1-r-s},$ and the coefficients are all
nonzero. Therefore there exists an $r$ such that  $\elt{u}{i}{r+s}\ot
\elt{v}{j}{\theta+d-1-r-s}\neq 0.$ In particular, $\elt{u}{i}{r+s}$
and $\elt{v}{j}{\theta+d-1-r-s}$ are nonzero, so $r+s\ppq \alpha-1$ and 
$\theta+d-1-r-s\ppq \beta-1.$ Therefore $(r+s)+(\theta+d-1-r-s)\ppq \alpha+\beta-2$.

On the other hand, $(r+s)+(\theta+d-1-r-s)=\theta+d-1.$ 

Now consider the case $\varpi<\theta<\alpha:$ then
$$\begin{array}{ll}
\alpha+\beta-2&= \varpi+d+1-2 \\&<\theta+d-1\\
&=
(r+s)+(\theta+d-1-r-s)\\&\ppq \alpha+\beta-2,
\end{array}$$
 and this gives a
contradiction, so $x_\theta$ must generate a simple module of
dimension $<d.$

Finally, assume that $\theta=\frac{\varpi}{2}$ with $\varpi$
even. Then either $\varpi>0$ and $2\theta=\varpi=\alpha+\beta-d-1,$ or
$\varpi=0$ and $\theta=0$ (since $\varpi\ppq 0$ and
$\theta=\frac{\varpi}{2}\pgq 0$).

Suppose  for a
contradiction  that $x_\theta$ generates a non-simple
module. Then the degree of $x_\theta$ must be ({see} Proposition \ref{projective})
$$\begin{array}{ll}
 -\rep{2(i+j)+(u+v)-1+2\theta}&=-\rep{(2i+u-1)+(2j+v-1)+1+\varpi}\\
&=-\rep{d-\alpha+d-\beta+1+\varpi}\\
&=-\rep{d}=-d,
\end{array}
$$
which is impossible since the smallest possible degree is $-d+1.$ 

Therefore $x_\theta$ also generates a simple module, and this module is
$L(u+v,i+j+\frac{\varpi}{2}),$ of dimension $d$ and so projective.

The sum of all the simple modules above is  contained in
$L(u,i)\ot L(v,j)$, and since the simple modules which occur are
non-isomorphic, the sum is direct.
\epf

We can now give the decomposition of $L(u,i)\ot L(v,j)$ when $\dim L(u,i) + \dim L(v,j) \ppq d+1:$ 

\bp\label{decompositionsmall} Assume that $\dim L(u,i) + \dim L(v,j) \ppq d+1.$ Then $$L(u,i)
\ot L(v,j) \cong \bigoplus _{\theta=0}^{\mathrm{min}\set{\dim L(u,i)-1
    ,\dim L(v,j)-1 }} L(u+v,i+j+\theta). $$ \ep

\bpf We know that   $L(u,i)
\ot L(v,j)$ contains  $\bigoplus _{\theta=0}^{\alpha-1} L(u+v,i+j+\theta). $ Computing the dimension of
this direct sum gives:  $$\dim \left(\bigoplus _{\theta=0}^{\alpha-1} L(u+v,i+j+\theta) \right)= \sum_{\theta=0}^{\alpha-1} \left(
d-\rep{-\alpha-\beta+1+2\theta}^-\right).$$ But we have $$-\alpha-\beta+1+2\theta\pgq
-(d+1)+1+2\theta \pgq -d $$ and $$-\alpha-\beta+1+2\theta\ppq
-\alpha-\beta+1+2(\alpha-1)=\alpha-\beta-1<0.$$ Therefore
$$\begin{array}{rcl}
\dim \left({\displaystyle\bigoplus _{\theta=0}^{\alpha-1} L(u+v,i+j+\theta)} \right)&=& {\displaystyle\sum_{\theta=0}^{\alpha-1}} \left(
d- (-\alpha-\beta+1+2\theta +d)\right)\\
&=& (\alpha+\beta-1)\alpha-2\frac{(\alpha-1)\alpha}{2}\\&=&\alpha\beta\\&=&\dim L(u,i)\ot L(v,j)
\end{array} $$
\epf


\subsubsection{Socle of $L(u,i)\ot L(v,j)$}

\bp \label{socle}  The socle of
$L(u,i)\ot L(v,j)$ is isomorphic to
$$\bigoplus_{\theta=\varsigma}^{\mathrm{min}\set{\dim L(u,i)-1;\dim
    L(v,j)-1}}L(u+v,i+j+\theta).$$ \ep

The proof follows from the next two lemmas; we first consider the case $0\ppq2\theta<\varpi:$

\bl Assume that $0\ppq2\theta<\varpi.$ Then $x_\theta$ does
not generate a simple module.\el

\bpfl  Note that we have $\varpi >0$ and therefore $\alpha+\beta>d+1.$
Suppose for a contradiction that $x_\theta$ does generate a
simple module. Then it would have dimension
$d-\rep{-\alpha-\beta+1+2\theta}^-=d-(-\alpha-\beta+1+2\theta+2d)=\alpha+\beta-d-1-2\theta.$
Therefore, the element $\gamma^{\alpha+\beta-d-1-2\theta}x_\theta$ would
be in the kernel of all arrows. 

But in the same way that we determined
the kernel of the action of $X$ ({see} Proposition \ref{kernel}), we can see that the intersection of the kernels of the actions of all the  arrows on $L(u,i)\ot L(v,j)$ has
dimension $\alpha$ and has a basis formed of elements $y_\pi$ for $\beta-1\ppq
\pi \ppq \alpha+\beta-2,$ with $\gr y_\pi=\pi$.

Therefore, we see that the element
$\gamma^{\alpha+\beta-d-1-2\theta}x_\theta$ must have a degree in
$L(u,i)\ot L(v,j)$ which is at least $\beta-1.$ However, $\gr
\gamma^{\alpha+\beta-d-1-2\theta}x_\theta
=\alpha+\beta-d-1-\theta=\beta-1-(d-\alpha+\theta)<\beta-1.$  Therefore $x_\theta$ cannot
generate a simple module.   \epfl

The remaining case is $\varpi<2\theta\ppq 2\varpi:$

\bl The socle of  the module generated by $x_\theta$  with $0\ppq2\theta<\varpi$ is a
simple module, generated by an element $x_{\theta'}$ with
$\varpi<2\theta'\ppq 2\varpi,$ and all the
remaining elements in the kernel of the action of $X$ can be obtained
in this way. \el

\bpfl The element $x_\theta$ is of the form $\tilde{D}_{u+v,p} $ for
some $p.$ To find $p,$ note  that $\gr x_\theta =i+j+\theta,$ so $e_{i+j+\theta}$
is the unique vertex acting as the identity on $\tilde{D}_{u+v,p}
$, and   therefore $i+j+\theta=p-\rep{2p+u+v-1}=\sigma_{u+v}(p)-d.$ So
we have $p=\sigma_{u+v}^2(p)-d=\sigma_{u+v}(d+i+j+\theta)-d.$

We can now compute the dimension of $\frac{\Gamma_u x_\theta}{\soc
  \Gamma_u x_\theta}$; this is equal to
$\rep{2p+u+v-1}=\varpi-2\theta$ using Proposition \ref{projective} and  the previous remarks.

Now consider $\gamma^{i+j+\theta}x_\theta$: we know that this
is an element in the kernel of the action of $X$.
Set $\theta':=\gr \gamma^{i+j+\theta}x_\theta=\varpi-\theta.$ Then
$2\theta'=\varpi+(\varpi-2\theta)>\varpi,$ so
$\gamma^{i+j+\theta}x_\theta=x_{\theta'}$ up to scalars, with $\varpi<2\theta'\ppq 2\varpi.$

The final claim follows from the fact that
$\#\set{\theta\,\mid\,0\ppq2\theta<\varpi }=\#\set{\theta'\,\mid\,\varpi<2\theta'\ppq 2\varpi}. $
\epfl


\subsubsection{Proof of the decomposition when $\alpha+\beta>d+1$}

We now know the socle of $L(u,i)\ot L(v,j),$ and we want to determine
the entire decomposition. For this we  determine a largest semisimple
summand which has no projective summands and prove that the remaining summands must be projective.

Consider the projective module $P(v,\sigma_{v}^{2p+1}(j)):$ $$\xymatrix@C=.3cm@=.2cm{&L(v,\sigma_v^{2p+1}(j))\ar@{-}[dl]\ar@{-}[dr]\\
 L(v,\sigma_v^{{2p}}(j))&&L(v,\sigma_v^{{2p}+2}(j))\\&L(v,\sigma_v^{2p+1}(j))\ar@{-}[ul]\ar@{-}[ur]}$$

We have an exact sequence $$0 \rightarrow L(v,
\sigma_v^{2p+1}(j)){\longrightarrow}  \Omega(L(v,
\sigma_v^{2p+1}(j)))\longrightarrow L(v,\sigma_v^{2p}(j))\oplus
L(v,\sigma_v^{2p+2}(j))\rightarrow 0.$$

Note that $\dim L(v,\sigma_v^{2p+1}(j))=d-\beta < \alpha$ since
$\alpha+\beta >d+1.$ Moreover, $\dim L(u,i)+\dim
L(v,\sigma_v^{2p+1}(j))=\alpha+d-\beta <d+1$ since $\alpha \ppq \beta.$
Therefore, by Proposition \ref{decompositionsmall},
taking the direct sum of these sequences over $p$ and tensoring on the
left by $L(u,i)$ gives:
$$
\begin{array}{ll}
0\rightarrow {\displaystyle\bigoplus_{p=0}^{\frac{n}{d}-1}
\bigoplus_{\theta=0}^{d-\beta-1} } L(u+v,i+\sigma_v^{2p+1}(j)+\theta)
\stackrel{f}{\longrightarrow} &
{\displaystyle\bigoplus_{p=0}^{\frac{n}{d}-1} 
\bigoplus_{\theta=0}^{d-\beta-1}} \Omega(L(u+v,i+\sigma_v^{2p+1}(j)+\theta)
) \oplus P  \\&\longrightarrow\left( {\displaystyle\bigoplus_{p=0}^{\frac{n}{d}-1}} L(u,i)\ot L(v,\sigma_v^{2p}(j))  \right)^2 
\rightarrow 0,
\end{array}
$$ where $P$ is a projective module.
 
Set $M:=\bigoplus_{p=0}^{\frac{n}{d}-1} L(u,i)\ot
L(v,\sigma_v^{2p}(j)).$ We know by Proposition \ref{socle} that the socle
of $M$ is
$\bigoplus_{p=0}^{\frac{n}{d}-1}\bigoplus_{\theta=\varsigma}^{d-\alpha-1}L(u+v,i+\sigma_v^{2p}(j)+\theta).$
It is straightforward to check that the summands in
$\bigoplus_{\theta=0}^{d-\beta-1}  L(u+v,i+\sigma_v^{2p+1}(j)+\theta)
$ and in $\soc (M)$ are pairwise non-isomorphic. 

We now determine the nonprojective part of $M.$ 

Take a  simple module $S:=L(u+v,i+\sigma_v^{2p+1}(j)+\theta)$ inside the
left-hand term. If $S$ is projective, then it must be a summand in
$P$ (since $\Omega(S)=0,$ or since $S$ is also injective), and does not occur in $M$. 

Now assume that $S$ is not projective. This simple module must embed via $f$ into either
$$\def\objectstyle{\scriptstyle}\xymatrix@=.2cm{L(u+v,\sigma_{u+v}^{-1}(i+\sigma_v^{2p+1}(j)+\theta))&&L(u+v,\sigma_{u+v}(i+\sigma_v^{2p+1}(j)+\theta))\\&S\ar@{-}[ur]\ar@{-}[ul]}$$
or an indecomposable projective summand $Q$ of $P$. Assume for a contradiction
that $S$ embeds into a projective $Q.$ Then the cokernel of this
embedding is
$$\def\objectstyle{\scriptstyle}\xymatrix@=.2cm{Q/S:&&&S\ar@{-}[dr]\ar@{-}[dl]\\&&L(u+v,\sigma_{u+v}^{-1}(i+\sigma_v^{2p+1}(j)+\theta))&&L(u+v,\sigma_{u+v}(i+\sigma_v^{2p+1}(j)+\theta))}.$$
This is therefore an indecomposable summand of $M^2,$ and hence of
$M. $ So the square of $Q/S$ must be a summand of $M^2.$ The second
copy of $Q/S$ must arise in the same way, so there are two copies of
the simple module $S$ in $\bigoplus_{\theta=0}^{d-\beta-1}  L(u+v,i+\sigma_v^{2p+1}(j)+\theta)
$, which contradicts the fact that the summands are pairwise
non-isomorphic.

Therefore $S$ embeds in
$$\def\objectstyle{\scriptstyle}\xymatrix@=.2cm{L(u+v,\sigma_{u+v}^{-1}(i+\sigma_v^{2p+1}(j)+\theta))&&L(u+v,\sigma_{u+v}(i+\sigma_v^{2p+1}(j)+\theta))\\&L(u+v,i+\sigma_v^{2p+1}(j)+\theta)\ar@{-}[ur]\ar@{-}[ul]}.$$
 So we have $M^2\cong P\oplus
 \bigoplus_{p=0}^{\frac{n}{d}-1}\bigoplus_{\theta=0}^{d-\beta-1}
 L(u+v,\sigma_{u+v}(i+\sigma_v^{2p+1}(j)+\theta))\oplus \bigoplus_{p=0}^{\frac{n}{d}-1}\bigoplus_{\theta=0}^{d-\beta-1}
 L(u+v,\sigma_{u+v}^{-1}(i+\sigma_v^{2p+1}(j)+\theta))$ which is
 isomorphic to  $\bigoplus_{p=0}^{\frac{n}{d}-1}\bigoplus_{\theta=0}^{d-\beta-1}
 L(u+v,\sigma_{u+v}(i+\sigma_v^{2p+1}(j)+\theta))^2$ (recall that
 $\sigma^2(\ell)=\ell+d$). 

We know the socle of $M,$ so we need to identify which  simple modules
in the socle we have obtained. We can see that
$\sigma_{u+v}(i+\sigma_v^{2p+1}(j)+\theta)=i+\sigma_v^{2p+2}(j)+(\alpha-1-\theta),$
and $\varpi+1=\alpha+\beta-d\ppq \pi:=\alpha-1-\theta \ppq \alpha-1.$ So $M\cong \bigoplus_{p=0}^{\frac{n}{d}-1}\bigoplus_{\pi=\varpi+1}^{\alpha-1}
 L(u+v,i+\sigma_v^{2p}(j)+\pi)\oplus P'$ with $P'$ projective. Since we know the socle of
 $M$, we get  $P'\cong\bigoplus_{p=0}^{\frac{n}{d}-1}
 \bigoplus_{\pi=\varsigma}^{\varpi}P(u+v,i+\sigma_v^{2p}(j)+\pi).$ 

Finally, since all the summands in $M$ are non-isomorphic, we conclude
using the socle of $L(u,i)\ot L(v,j).$

\subsection{Tensor product of modules of odd length}\label{tensorodd}

We know from  the description of the modules of odd length \ref{desc}{\bf (\ref{descodd})} that all the modules of odd
length are syzygies of simple modules. Moreover,  we
have $$
\begin{array}{ll}
\Omega^k(M)\ot \Omega^\ell(N)&\cong \Omega^k(M\ot \Omega^\ell(N))\oplus
\mbox{ (projectives) } \\&\cong \Omega^k(\Omega^\ell(N)\ot M)\oplus
\mbox{ (projectives) }\\&\cong
\Omega^k(\Omega^\ell(N\ot M))\oplus
\mbox{ (projectives) }\\&\cong \Omega^{k+\ell}(M\ot N)\oplus
\mbox{ (projectives)}.
\end{array}
$$ Therefore, from the results above applied with $M$ and $N$ simple, we can decompose the tensor
product of any modules of odd length into a direct sum of modules of
odd length, up to projectives.

\subsection{Tensor product of a string module of even length with a module of
  odd length}

In this section, we study the tensor product of a string  module of even length with a module of
  odd length, and by \ref{desc}{\bf (\ref{descodd})} and Section \ref{tensorodd}   it is enough to consider the
  tensor product of a  string module of even length with a simple module.
To simplify notations, let $\frak{I}$ denote the set $\set{0\ppq \theta \ppq {\rm min}\set{\dim L(u,i),\dim
    L(v,j)}-1}$ if $\dim L(u,i)+\dim
    L(v,j)<d+1$, and the set $\set{\varpi +1\ppq \theta \ppq {\rm min}\set{\dim L(u,i),\dim
    L(v,j)}-1}$ if $\varpi:=\dim L(u,i)+\dim
    L(v,j)-(d+1)\pgq 0.$


\bp\label{twosimple} The tensor product of a string module of length
two with a simple module decomposes as follows: $$M_2^+(u,i)\ot L(v,j) \cong 
\displaystyle{\bigoplus_{\theta\in\frak{I}}M_2^+(u+v,i+j+\theta)\oplus {\rm projective}
}$$and  $$M_2^-(u,i)\ot L(v,j) \cong 
\displaystyle{\bigoplus_{\theta\in\frak{I}}M_2^-(u+v,i+j+\theta)\oplus {\rm projective}}.$$ 
\ep

\bpf
We prove it for $M_2^+(u,i).$ Consider the exact sequence $$0\rightarrow L(u,\sigma_u(i))\rightarrow
M_2^+(u,i) \rightarrow L(u,i)\rightarrow 0$$ and tensor it with $L(v,j).$ 
We know the decomposition of the two outside terms from Theorem
\ref{tensorsimples}, so $M_2^+(u,i)\ot L(v,j)=W\oplus P$ where $P$ is a
projective module and the composition factors of $W$ are the non-projective
simple summands in $L(u,i)\ot L(v,j)$ and $L(u,\sigma_u(i))\ot L(v,j),$ which
are known.

We now show that the non-projective  summands in $M_2^+(u,i)\ot L(v,j)$ have
length 2. We know that $M_2^+(u,i)$ is periodic (from \ref{desc}{\bf (\ref{descstring})}), so $\Omega^m(M_2^+(u,i))\cong M_2^+(u,i)$ for some $m \in \mathbb{Z}.$
Then, tensoring a minimal projective resolution of $M_2^+(u,i)$ by $L(v,j)$ and
using the fact that $P\ot L(v,j)$ is projective if $P$ is projective
\cite[Proposition 2.1]{GMS} gives  $\Omega^m(M_2^+(u,i)\ot L(v,j))\cong M_2^+(u,i)\ot L(v,j)\oplus
{\rm projective}.$ So any summand of $M_2^+(u,i)\ot L(v,j)$ is either
 periodic, and therefore of even length, or projective. From the
 description of the indecomposable modules of even length
 (\ref{desc}), we know that all the simple summands in the top have
 the same dimension $N$ and all the simple summands in the socle have
 the same dimension $d-N.$ 
It is easy to check that the (non-projective) simple modules which are
summands in $L(u,i)\ot L(v,j)$ all have different dimensions, and that
for each dimension $N$ that occurs there is exactly one simple summand
in $L(u,\sigma_u(i))\ot L(v,j)$ which has dimension $d-N,$ and
moreover  that the
simple modules pair off to give indecomposable summands of length 2,
whose top and socle are  fully determined by the decompositions of $L(u,i)\ot L(v,j)$
and $L(u,\sigma_u(i))\ot L(v,j)$.

We now need to decide whether these summands of length 2 are string or
band modules. If $n\neq d,$ then there are no band modules of length 2 so
it is clear. 

Now assume that $n=d,$ and let us write temporarily $C_0^{1+}(u,i)$
for $M_2^+(u,i).$ For any $\lambda\in k,$ the module
$C_\lambda^{1+}(u,i)$ is characterised by the following pullback diagram:
$$\xymatrix{0\ar[r]&L(u,\sigma_u(i))\ar[r]\ar@{=}[d]&C_\lambda^{1+}(u,i)\ar[r]\ar[d] &L(u,i)\ar[d]^{\varphi_\lambda}\ar[r]&0\\
            0\ar[r]&L(u,\sigma_u(i))\ar[r]
            &P(u,\sigma_u(i))\ar[r]&\Omega^{-1}(L(u,\sigma_u(i)))\ar[r] &0}$$ with $\varphi_\lambda(\tilde{H}_{u,i})=\lambda \overline{b}_0D_{u,\sigma_u(i)}+{b}_1D_{u,\sigma_u(i)}$. 

The summands in $C_0^{1+}(u,i)\ot L(v,j)$ are of the form
$C_\lambda^{1+}(u+v,i+j+\theta).$ Consider such a summand: then $\varphi_\lambda:L(u+v,i+j+\theta)\rightarrow
\Omega(L(u+v,\sigma_{u+v}(i+j+\theta)))$ is the restriction of
$\varphi_0\ot id: L(u,i)\ot L(v,j)\rightarrow \Omega(L(u,i))\ot
L(v,j)$ to the summand $L(u+v,i+j+\theta).$ Applying $b_0$ to the
first term in $\varphi_\lambda(\tilde{H}_{u+v,i+j+\theta})$ gives
$\lambda b_0\bar{b}_0 D_{u+v,i+j+\theta}$, and applying $\bar{b}_1$ to
the second term gives $\bar{b}_1b_1 D_{u+v,i+j+\theta}= b_0\bar{b}_0
D_{u+v,i+j+\theta}.$ So the first one is $\lambda$ times the
second, and they are both nonzero. Now apply the same procedure with $\varphi_0\ot \id:$ the
first term will be 0, and the second non-zero. Hence we must have
$\lambda=0$ and therefore the summands in $M_2^+(u,i)\ot L(v,j)$ which
are of the form $C_\lambda^{1+}(u+v,i+j+\theta)$ are in fact
of the form $M_2^+(u+v,i+j+\theta).$

We then consider all the cases: 
\begin{enumerate}[$\bullet$]
\item $\dim L(u,i)+\dim L(v,j)<  d+1$ and $\dim L(u,i) \ppq \dim L(v,j);$
\item $\dim L(u,i)+\dim L(v,j) < d+1$ and $\dim L(u,i) > \dim L(v,j);$
\item $\dim L(u,i)+\dim L(v,j) \pgq d+1$ and $\dim L(u,i) > \dim L(v,j);$
\item $\dim L(u,i)+\dim L(v,j) \pgq  d+1$ and $\dim L(u,i) \ppq \dim L(v,j)$
\end{enumerate} in order to determine the exact bounds (for $\theta$) of the decomposition up to projectives.

\epf

\bt\label{stringsimple} The tensor product of a string module of even
length with  a simple module decomposes as follows: $$M_{2\ell}^+(u,i)\ot L(v,j) \cong 
\displaystyle{\bigoplus_{\theta\in\frak{I}}M_{2\ell}^+(u+v,i+j+\theta)\oplus
  {\rm projective}}  
$$ and  $$M_{2\ell}^-(u,i)\ot L(v,j) \cong 
\displaystyle{\bigoplus_{\theta\in\frak{I}}M_{2\ell}^-(u+v,i+j+\theta)\oplus {\rm projective}}
$$ 
\et

\bpf We work by induction on $\ell,$ using Proposition \ref{arandstm}
(\ref{arstm}), which we apply here with 
 $M=M^+_{2\ell}(u,i+d)$ and $N=L(v,j)$. We know that $N$ is a
 splitting trace module (Proposition \ref{stm}) and we have the
 Auslander-Reiten sequence $$\A(M^+_{2\ell}(u,i+d)):\ \ 
0\rightarrow M^+_{2\ell}(u,i)\rightarrow 
\begin{array}{c}
M^+_{2\ell+2}(u,i)\\\oplus\\
M^+_{2\ell-2}(u,i+d)
\end{array}\rightarrow M^+_{2\ell}(u,i+d)\rightarrow 0
$$ (where $M^+_0(u,i):=0$).  The summands in $M^+_{2\ell+2}(u,i+d)\ot
L(v,j)$ are not splitting trace modules (they have even length), therefore, by Proposition \ref{arandstm}
(\ref{arstm}), the sequence $\A(M^+_{2\ell}(u,i+d))\ot L(v,j)$ is
split exact. We know the decomposition of the
tensor product of three of the terms with $L(v,j)$, so induction
yields the fourth.
\epf

\br We can now find the tensor product of any  string module of even
length with any string module of odd length up to
projectives,
using the results and methods of the previous sections.
\er

\subsection{Tensor product of string modules of even length}


\bp The tensor product of two string modules of length two of the same
type decomposes as follows: $$M_2^+(u,i)\ot M_2^+(v,j) \cong 
\displaystyle{\bigoplus_{\theta\in\frak{I}}[M_2^+(u+v,i+j+\theta)\oplus
M_2^+(u+v,\sigma_{u+v}(i+j+\theta))] \oplus {\rm projective}}$$ and  $$M_2^-(u,i)\ot M_2^-(v,j) \cong 
\displaystyle{\bigoplus_{\theta\in\frak{I}}[M_2^-(u+v,i+j+\theta)\oplus
M_2^-(u+v,\sigma_{u+v}^{-1}(i+j+\theta))] \oplus {\rm projective}}.$$  
 \ep

\bpf
We start with the exact sequence $0 \rightarrow
L(v,\sigma_v^{-1}(j))\rightarrow\Omega^{-1}(L(v,j))\rightarrow
M_2^+(v,j)\rightarrow 0$.  Tensoring on the left by $M_2^+(u,i)$ gives
again an exact sequence, and using Proposition \ref{twosimple} we get the
exact sequence: 
$$0\rightarrow \oplus_\rho
M_2^+(u+v,i+\sigma_v^{-1}(j+\rho))\rightarrow
\Omega^{-1}(\oplus_\theta M_2^+(u+v,i+j+\theta))\oplus P\rightarrow
M_2^+(u,i)\oplus M_2^+(v,j)\rightarrow 0$$ where $P$ is a projective module.
Note that the middle term is equal to $\oplus_\theta
M_2^+(u+v,\sigma_{u+v}(i+j+\theta))\oplus P$ and that
$i+\sigma_v^{-1}(j)+\rho=\sigma_{u+v}^{-1}(i+j+\pi)$ with $\pi$ in the
same range as $\theta$ (for the latter, we consider all the four cases
that occurred in the proof of Proposition \ref{twosimple}). So we
have $$0\rightarrow \oplus _\pi
M_2^+(u+v,\sigma_{u+v}^{-1}(i+j+\pi))\rightarrow \oplus_\theta
M_2^+(u+v,\sigma_{u+v}(i+j+\theta))\oplus P\rightarrow
M_2^+(u,i)\oplus M_2^+(v,j)\rightarrow 0.$$

We now prove that $M_2^+(u+v,\sigma_{u+v}^{-1}(i+j+\pi))$ embeds into
the projective $P:$ if not, then the socle of
$M_2^+(u+v,\sigma_{u+v}^{-1}(i+j+\pi))$, which is $L(u+v,i+j+\pi)$,
embeds into one of the modules of length two, so we would have
$L(u+v,i+j+\pi)=L(u+v,\sigma^2_{u+v}(i+j+\theta))=L(u+v,i+j+\theta+d)$
and hence $\pi=\theta+d,$ which is impossible since $\pi$ and $\theta$
are in the same range, which has at most $d-1$ elements.

Therefore $M_2^+(u+v,\sigma_{u+v}^{-1}(i+j+\pi))$ embeds into
the projective $P$ and gives $M_2^+(u+v,i+j+\pi)$ in the
quotient. This means that in the quotient we have a projective, $M_2^+(u+v,i+j+\pi)$
and $M_2^+(u+v,\sigma_{u+v}(i+j+\theta))$ with $\theta$ and $\pi$ in
the same range. We cannot have $i+j+\pi=\sigma_{u+v}(i+j+\theta)$
(again, consider cases as in the proof of Proposition \ref{twosimple}) so
these modules of length two do not link up to give modules of length
four. They must therefore occur as summands.

\epf

\bp The module $M_2^-(u,i)\ot M^+_2(v,j)$ is projective.
\ep

\bpf If we tensor the exact sequence $0\rightarrow
L(v,\sigma_v^{-1}(j))\rightarrow \Omega^{-1}(L(v,j))\rightarrow
M_2^+(v,j)\rightarrow 0$ on the left by $M_2^-(u,i)$, we get the exact
sequence $0\rightarrow\oplus_\pi
M_2^-(u+v,\sigma_{u+v}(i+j+\pi))\rightarrow\oplus_\theta
M_2^-(u+v,\sigma_{u+v}^{-1}(i+j+\theta))\oplus P\rightarrow
M_2^-(u,i)\ot M^+_2(v,j)\rightarrow 0$ with $\pi$ and $\theta$ in the
same range. So $M_2^-(u,i)\ot M^+_2(v,j)$ must be a sum of modules of the form $M_2^-(u+v,\ell)$ and a
projective.

A similar argument starting with the exact sequence $0\rightarrow
L(u,\sigma_u(i))\rightarrow \Omega^{-1}(L(u,i))\rightarrow
M_2^-(u,i)\rightarrow 0$ implies that  $M_2^-(u,i)\ot M^+_2(v,j)$ must be a sum of modules of the form $M_2^+(u+v,\ell)$ and a
projective.

Combining these results implies that $M_2^-(u,i)\ot M^+_2(v,j)$ is projective since the decomposition into indecomposable summands is unique and we can distinguish modules of the form $M_2^+(w,\ell)$ and modules of the form $M_2^-(w,\ell)$ by viewing them as $\land$-modules and checking whether as such they are projective or not.

\epf

\bt If $\ell$ and $t$ are two positive integers, then
$M_{2\ell}^+(u,i)\ot M_{2t}^-(v,j)$ is projective, and we have the
decompositions: $$M_{2\ell}^+(u,i)\ot M_{2t}^+(v,j) \cong 
\displaystyle{\bigoplus_{\theta\in\frak{I}}\bigoplus_{p=0}^{2\ell-1}\bigoplus_{r=0}^{t-1}}M_2^+(u+v,\sigma_{u+v}^{p+2r}(i+j+\theta))
\oplus {\rm projective}$$ and $$M_{2\ell}^-(u,i)\ot M_{2t}^-(v,j) \cong 
\displaystyle{\bigoplus_{\theta\in\frak{I}}\bigoplus_{p=0}^{2\ell-1}\bigoplus_{r=0}^{t-1}}M_2^-(u+v,\sigma_{u+v}^{-(p+2r)}(i+j+\theta)) \oplus {\rm projective}$$

\et

\bpf
We work by induction on $\ell$ and $t,$  using Proposition \ref{arandstm}
(\ref{arnotstm}), since we know that modules of even length are not
splitting trace modules (Proposition \ref{notstm}), so we will obtain
split exact sequences.

Let us consider for instance $M_{2\ell}^+(u,i)\ot M_{2t}^+(v,j)$. We first set $t=1$ and consider the Auslander-Reiten sequence for
$M_{2\ell}^+(u,i+d).$ An induction on $\ell$ gives
$M_{2\ell+2}^+(u,i+d)\ot M_2^+(v,j)$ and hence
$M_{2\ell}^+(u,i+d)\ot M_2^+(v,j)$ for all $\ell.$ 

We can then do an induction on $t,$ considering the Auslander-Reiten
sequence for $M_{2t}^+(v,j+d)$, to get the result.
\epf

\subsection{Tensor products of band modules with modules of odd
  length}

Consider the tensor product $C_\lambda^{1+}(u,i)\ot L(v,j)$. Using a
method similar to the proof of Proposition \ref{twosimple}, we can see
that this tensor product decomposes as $\bigoplus_\theta
C_{\mu_\theta}^{1+}(u,i).$ We must then determine the
parameters $\mu_\theta$. We will describe the method on an example,
which will give an algorithm to determine these parameters for any
given example, and illustrate the complexity of the general case.

\be\label{examplebandsimple} Let us consider the tensor product $C^{1+}_\lambda(1,5)\ot L(0,2)$
when $n=d=6.$ As we have said above, we know that this tensor product
decomposes as $C^{1+}_{\mu_0}(1,1)\oplus
C^{1+}_{\mu_1}(1,2)\oplus L(1,3)$ for some $\mu_0$ and $\mu_1$ to be
determined.

The blocks involved are 
 $$\xymatrix@1{ L(1,5) \ar@/^2pc/@<+1ex>[rrr]^{b_0=\gamma^2}
  \ar@/^1pc/@<+0.5ex>[rrr]^{\overline{b}_1=X^4} &&& L(1,1) \ar@/^1pc/@<+0.5ex>[lll]^{b_1=\gamma^4}
  \ar@/^2pc/@<+1ex>[lll]^{\overline{b}_0=12X^2} 
&\mbox{ and }& L(1,2) \ar@/^2pc/@<+1ex>[rrr]^{b_0'=\gamma^2}
  \ar@/^1pc/@<+0.5ex>[rrr]^{\overline{b}_1'=X^4} &&& L(1,4) \ar@/^1pc/@<+0.5ex>[lll]^{b_1'=\gamma^4}
  \ar@/^2pc/@<+1ex>[lll]^{\overline{b}_0'=12X^2} }.$$

\paragraph{Determine $\mu_0:$} Recall that,  for any $\lambda\in k,$ the module
$C_\lambda^{1+}(u,i)$ is characterised by the following pullback diagram:
$$\xymatrix{0\ar[r]&L(u,\sigma_u(i))\ar[r]\ar@{=}[d]&C_\lambda^{1+}(u,i)\ar[r]\ar[d] &L(u,i)\ar[d]^{\varphi_\lambda}\ar[r]&0\\
            0\ar[r]&L(u,\sigma_u(i))\ar[r]
            &P(u,\sigma_u(i))\ar[r]&\Omega^{-1}(L(u,\sigma_u(i)))\ar[r] &0}$$ where the homomorphism $\varphi_\lambda$ is determined by $\varphi_\lambda(\tilde{H}_{u,i})=\lambda \overline{b}_0D_{u,\sigma_u(i)}+{b}_1D_{u,\sigma_u(i)}$. 

We first find $Z_0$ in $L(1,5)\ot L(0,2)$ corresponding to
  $\tilde{H}_{11}$ in $L(1,1)$ (up to scalar multiples): this satisfies $e_1Z_0=Z_0$ so we must
  have $Z_0=\tilde{H}_{15}\ot \tilde{H}_{02}$, and we check that
  $XZ_0=0.$ 

We then apply $\varphi_\lambda\ot \id$ to $Z_0:$ this is equal to
  $(x_0,y_0)$ with $x_0=\lambda (\bar{b}_0 D_{11})\ot \tilde{H}_{02}=12
  \lambda \tilde{D}_{11}\ot \tilde{H}_{02}$ and $y_0=(b_1D_{11})\ot
  \tilde{H}_{02}=-12 H_{11}\ot\tilde{H}_{02}. $ 

We now need to lift $x_0$ and $y_0$ to $P(1,1)\ot L(0,2):$
\begin{enumerate}[$\bullet$]
\item Since $Xx_0=0,$ we know that $x_0$ corresponds to $\tilde{D}_{15},$
  so we do not change $x_0.$
\item We look for $\bar{y}_0=-12 H_{11}\ot\tilde{H}_{02}+\alpha X
  H_{11}\ot a\tilde{H}_{02} +\beta X^2H_{11}\ot \gamma^2\tilde{H}_{02}
  $. The element $\bar{y}_0$ is characterised by $X\bar{y}_0\neq 0$ and
  $aX\bar{y}_0=0$ (it corresponds to $H_{15}$). This
  gives the equations $-12q^2+\alpha(2q+q^2)+2\beta(1+q)=0$ and
  $2(q^2-1)\beta-\alpha =0,$ so finally $$\bar{y}_0=-12H_{11}\ot\tilde{H}_{02}+4(1+q) X
  H_{11}\ot a\tilde{H}_{02} -2q X^2H_{11}\ot \gamma^2\tilde{H}_{02}.$$
\end{enumerate}

We know that
  $\varphi_{\mu_0}(\tilde{H}_{11})=(\mu_0\bar{b}_1D_{15},b_0D_{15})$
  corresponds to $(x_0,\bar{y}_0),$ so $b_1x_0=\mu_0 \bar{b}_0\bar{y}_0.$

Calculating $b_1x_0$ gives $12\lambda[6X^2H_{11}\ot
  \tilde{H}_{02}-2(q+q^2)X^3H_{11}\ot a\tilde{H}_{02}-2
  \tilde{H}_{11}\ot \gamma^2 \ot \tilde{H}_{02}]$, and calculating
  $\bar{b}_0\bar{y}_0$ gives $12[-6 X^2H_{11}\ot
  \tilde{H}_{02}+2(q+q^2)X^3H_{11}\ot a\tilde{H}_{02}+2
  \tilde{H}_{11}\ot \gamma^2 \ot \tilde{H}_{02}]$, so finally
  $\mu_0=-\lambda. $

\paragraph{Determine $\mu_1:$} We proceed in the same way:
$\tilde{H}_{12}$ corresponds to $Z_1=(1+q)a\tilde{H}_{15}\ot \tilde{H}_{02}+\tilde{H}_{15}\ot
  a\tilde{H}_{02}$ in $L(1,5)\ot L(0,2)$, and applying
  $\varphi_\lambda\ot \id$ gives $(x_1,y_1)$  where
  $x_1=12\lambda[(q+q^2)XD_{11}\ot \tilde{H}_{02}+\tilde{D}_{11}\ot
  a\tilde{H}_{02}]$ and $y_1=12[-(q+q^2)F_{11}\ot \tilde{H}_{02}-H_{11}\ot
  a\tilde{H}_{02}]$ (note that
  $\varphi_\lambda(a\tilde{H}_{15})=a\varphi_\lambda(\tilde{H}_{15})$). Lifting to $P(1,1)\ot L(0,2)$ does not change $x_1$, and $y_1$ becomes $\bar{y}_1=12[-(q+q^2)F_{11}\ot \tilde{H}_{02}-H_{11}\ot
  a\tilde{H}_{02}+XH_{11}\ot \gamma^2\tilde{H}_{02}].$

Then, since $\varphi_{\mu_1}(\tilde{H}_{12})=(\mu_1\bar{b}'_0D_{14},b'_1D_{14})$
  corresponds  to a scalar multiple of $(x_1,\bar{y}_1)$, we have
  $b'_0x_1=\mu_1\bar{b}'_1\bar{y}_1,$ and calculating each of the terms
  in this identity yields $\mu_1=\lambda$.
\ 

Therefore $C_\lambda^{1+}(u,i)\ot L(v,j)=C^{1+}_{-\lambda}(1,1)\oplus
C^{1+}_{\lambda}(1,2)\oplus L(1,3)$.

\ee

Using the method described in this example, we can calculate the
tensor product of a band module of minimal length $\frac{2n}{d}$ with
a simple module, and hence using Auslander-Reiten sequences and the
fact that modules of odd length are syzygies of simple modules, as in
the previous sections, we can determine the tensor product of any band
module with a module of odd length.

\subsection{Tensor products of band modules with other modules of even
  length}

Here again, we will describe the method on an example:

\be We shall consider the tensor product $C_\lambda^{1+}(1,5)\ot
C_\mu^{1+}(0,2)$ with $n=d=6.$  Since the modules we are tensoring
are periodic, their summands must also be periodic or projective, and
hence of even length. From Example \ref{examplebandsimple} and a
similar calculation which gives $C_\lambda^{1+}(1,5)\ot
L(0,5)=C_{-\lambda}^{1+}(1,4)\oplus C_\lambda^{1+}(1,5)\oplus L(1,0),$ we
see that the non-projective composition factors of $C_\lambda^{1+}(1,5)\ot
C_\mu^{1+}(0,2)$ are $L(1,1),$ $L(1,1)$, $L(1,5)$, $L(1,5)$,
$L(1,2)$, $L(1,2)$,  $L(1,4)$  and $L(1,4)$. Note that the first four
composition factors and the last four composition factors must be in
separate summands of  $C_\lambda^{1+}(1,5)\ot
C_\mu^{1+}(0,2)$, since $\sigma_1(1)=5$ and $\sigma_1(5)=1$,
neither of which is equal to 2 or 4. Call these summands $M$ and $N.$

\paragraph{Determine the top of  $C_\lambda^{1+}(1,5)\ot
C_\mu^{1+}(0,2)$:}  We shall use here the following result from
\cite[Propositions 1.1 and 1.2]{GMS}: $$\Hom(A\ot B,C)\cong
\Hom(A,B^*\ot C)$$ for any $\dland$-modules $A$, $B$ and $C,$ where
$B^*$ is the $k$-dual of $B.$
Moreover, we can see that
$C_\lambda^{1+}(u,i)^*\cong C_\lambda^{1+}(u,\sigma_u(i))$ (note that
it is easy to see, using Proposition \ref{simples}, that
$L(u,i)^*\cong L(1-u,\sigma_u(i))$). 

The module $L(1,1)$ is in the top of  $C_\lambda^{1+}(1,5)\ot
C_\mu^{1+}(0,2)$ if, and only if,  $\Hom(C_\lambda^{1+}(1,5)\ot
C_\mu^{1+}(0,2),L(1,1))$ is nonzero. We have $$
\begin{array}{rcl}
\Hom(C_\lambda^{1+}(1,5)\ot
C_\mu^{1+}(0,2),L(1,1))&\cong& \Hom(C_\lambda^{1+}(1,5),
C_\mu^{1+}(0,5)\ot L(1,1))\\
&\cong& \Hom(C_\lambda^{1+}(1,5),C_{\nu_1}^{1+}(1,1))\oplus
\Hom(C_\lambda^{1+}(1,5),C_{\nu_2}^{1+}(1,2))\\&&\  \oplus
\Hom(C_\lambda^{1+}(1,5),L(1,3))\\
&\cong& k\oplus 0\oplus 0=k.
\end{array}
$$ Hence $L(1,1)$ occurs once in the top.

Similarly, $L(1,2)$ occurs once in the top.
Now $$
\begin{array}{rcl}
\Hom(C_\lambda^{1+}(1,5)\ot
C_\mu^{1+}(0,2),L(1,5))&\cong& \Hom(C_\lambda^{1+}(1,5),
C_\mu^{1+}(0,5)\ot L(1,5))\\
&\cong& \Hom(C_\lambda^{1+}(1,5),C_{\nu_3}^{1+}(1,4))\oplus
\Hom(C_\lambda^{1+}(1,5),C_{\nu_4}^{1+}(1,5)) \\&&\  \oplus
\Hom(C_\lambda^{1+}(1,5), L(1,0))\\
&\cong & \Hom(C_\lambda^{1+}(1,5),C_{\nu_4}^{1+}(1,5))\\
&=& \begin{cases}k &\mbox{ if $\nu_4=\lambda$}\\0 & \mbox{ otherwise}, \end{cases}
\end{array}
$$ since any homomorphism must be an isomorphism. So we need to determine $\nu_4.$ Proceeding as in Example
\ref{examplebandsimple} gives $\nu_4=-\dfrac{1+q}{2}\mu.$ Hence
$L(1,5)$ occurs in the top if, and only if,
$\lambda=-\dfrac{1+q}{2}\mu.$

Similarly, $L(1,4)$   occurs in the top if, and only if,
$\lambda=\dfrac{1+q}{2}\mu.$

\paragraph{Determine the socle of $C_\lambda^{1+}(1,5)\ot
C_\mu^{1+}(0,2)$:} This is similar, considering $\Hom(L(u,i),C_\lambda^{1+}(1,5)\ot
C_\mu^{1+}(0,2))$ for each of the simples $L(u,i)$ involved. We
find that $L(1,1)$ and $L(1,2)$ both occur once in the socle, $L(1,5)$
occurs in the socle if, and only if, $\lambda=-\dfrac{1+q}{2}\mu,$ and
$L(1,4)$ occurs in the socle if, and only if,
$\lambda=\dfrac{1+q}{2}\mu.$

\paragraph{First case: $\lambda\neq \pm\dfrac{1+q}{2}\mu$:} Then
$L(1,4)$ and $L(1,5)$ do not occur in the top or the socle of $C_\lambda^{1+}(1,5)\ot
C_\mu^{1+}(0,2)$. Therefore the two summands we considered above
each have a
simple top and a simple socle, so  $C_\lambda^{1+}(1,5)\ot
C_\mu^{1+}(0,2)=P(1,1)\oplus P(1,2)\oplus {\rm projective}$ is a
projective module.

\paragraph{Second case:  $\lambda=\dfrac{1+q}{2}\mu$:} Then the
summand $M$ is as in the
first case and must be $P(1,1),$ but the summand $N$ differs: $L(1,2)
$ and $L(1,4)$ both  occur
once in the top and once in the socle, and since the summands must
have even length, $N$ must be $C^{1+}_{\nu_5}(1,2)\oplus C^{1+}_{\nu_6}(1,4).$

In order to determine $\nu_5$ and $\nu_6,$ consider the following
exact sequence: $$0\rightarrow L(0,5)\rightarrow C^{1+}_\mu(0,2)\rightarrow
L(0,2)\rightarrow 0.$$ Tensoring with $C^{1+}_\lambda(1,5)$ gives the
exact sequence: $$0\rightarrow C^{1+}_{-\lambda}(1,4)\oplus
C^{1+}_{\lambda}(1,5)\oplus L(1,0)\rightarrow C_\lambda^{1+}(1,5)\ot
C_\mu^{1+}(0,2)\rightarrow
C^{1+}_{-\lambda}(1,1)\oplus C^{1+}_{\lambda}(1,2)\oplus L(1,3)\rightarrow
0,$$ so we know, by considering the composition factors of $C^{1+}_{-\lambda}(1,4)$ and of $C_\lambda^{1+}(1,5)\ot
C_\mu^{1+}(0,2)$,  that $ C^{1+}_{-\lambda}(1,4)$ embeds into $C^{1+}_{\nu_5}(1,2)\oplus C^{1+}_{\nu_6}(1,4)$. So we have an embedding $\iota:C^{1+}_{-\lambda}(1,4) \rightarrow C^{1+}_{\nu_5}(1,2)\oplus C^{1+}_{\nu_6}(1,4) $ and a projection $\pi: C^{1+}_{\nu_5}(1,2)\oplus C^{1+}_{\nu_6}(1,4) \rightarrow C^{1+}_{\nu_6}(1,4) $. The composition $\pi\iota$ is non-zero on the socle $L(1,2),$ hence injective, so it must be an isomorphism. Therefore we have $\nu_6=-\lambda.$ A similar argument, using the fact that $C^{1+}_{\lambda}(1,2) $ is a  quotient of   $C_\lambda^{1+}(1,5)\ot
C_\mu^{1+}(0,2)$, shows that $\nu_5=\lambda.$

Finally $C_\lambda^{1+}(1,5)\ot
C_\mu^{1+}(0,2)= C^{1+}_{\lambda}(1,2)\oplus
C^{1+}_{-\lambda}(1,4)\oplus {\rm projective}.$

\paragraph{Third case:  $\lambda=-\dfrac{1+q}{2}\mu:$} This is similar
to the second case, and we get $C_\lambda^{1+}(1,5)\ot
C_\mu^{1+}(0,2)=C^{1+}_{-\lambda}(1,1)\oplus C^{1+}_{\lambda}(1,5)\oplus {\rm projective}.$

\ee

Other situations (such as $d\neq n$ or tensor products of different
types of band modules of minimal length) are treated similarly. Moreover, the
case of the tensor product of a band module with a string module of
 even length can also be treated in this way, setting the parameter to
 zero for the string module. Finally,
the general case (that is, when the lengths are not minimal) can be
deduced by induction, using Auslander-Reiten sequences.


\section{Interpretation of simple and projective $\dland$-modules as Hopf
  bimodules over $\land$}

In this section, we describe some Hopf bimodules over $\land:$ given
any finite-dimensional Hopf algebra $H,$ there
is an equivalence of categories between modules over the Drinfel'd
double $\D(H)$ and Hopf bimodules over $H$ ({see}
\cite{kassel,montgomery,rosso} for instance). Recall:

\bd Let $H$ be a Hopf algebra. A {\bf Hopf bimodule} over $H$ is an
$H$-bimodule $M$ which is an $H$-bicomodule such that the comodule
structure maps $M\rightarrow H\ot M$ and $M\rightarrow M\ot H$ are
$H$-bimodule homomorphisms ($H\ot M$ and $M\ot H$ are endowed with
diagonal $H$-bimodule structures).  We use Sweedler's notation for the comultiplication of $H,$ so $\Delta(h)=h^{(1)}\ot h^{(2)}.$ \ed

The functor from $\D(H)$-modules to Hopf bimodules over $H$ can be
described as follows: let $V$ be a module over $\D(H).$ Consider the
vector space $M:=H\ot V;$ as a right $H$-module, it is free (of rank
$\dim V$), that is, for any $h,a\in H$ and $v\in V,$ we have  $(h\ot v)\cdot a=ha\ot v.$ As a left $H$-comodule,
it is free, that is, $h\ot v \mapsto h^{(1)}\ot (h^{(2)}\ot v). $
The left module structure on $M$ is diagonal, that is, $a\cdot (h\ot
v)=a^{(1)}h\ot a^{(2)}v.$

Finally, to describe the right comodule structure, we need to fix a
basis $\set{b_i}$ of $H$ and $\set{b_i^*}$ its dual basis; then
$\delta_R(h\ot v)=\sum_i h^{(1)}\ot b_i^* \cdot v \ot b_i\cdot h^{(2)}.$
Note that $V$ is a left $H^*$-module and a left $H$-module (restricting the action of
$\D(H)$ to $H^{*cop}$ and to $H$). 

If we take $H=\land,$ we have a basis given by
$\set{q^{-im}\frac{1}{m!_{\! q}}\gamma_i^m\,\mid\, i\in \zz_n, 0\ppq
  m\ppq d-1}.$ We refer the reader to the notation and correspondence
in Section 2
(page \pageref{taft}) in order to see that the dual basis is
 $\set{X^mG^i\,\mid \, i\in \zz_n, 0\ppq m\ppq
  d-1}.$ Hence, if $V$ is a module over $\dland,$ we have $$\delta_R(h\ot v)=\sum_{\tiny 
\begin{array}{c}i\in \zz_n \\ 0\ppq m \ppq d-1
\end{array}}
(h^{(1)}\ot X^mG^{i}v )\ot q^{-im}\frac{1}{m!_{\! q}}\gamma_i^mh^{(2)}.$$

Therefore, if $V$ is a module over $\dland,$ the corresponding Hopf
bimodule over $\land$ is free as a right module and left comodule, so
as such it is isomorphic to $\land ^{\dim V}$, and we now need to
describe the left module and right comodule structures on $\land ^{\dim
  V}$ using the definitions above.

\subsection{Hopf bimodule corresponding to $L(u,j)$}

Let $N$ be the dimension of $L(u,j).$ Let us first describe the left
$\land$-module structure of $\land \ot L(u,j).$ 

Define $\tau_{k}^{(p)}(\gamma_i^m)=\begin{pmatrix} m\\
  p\end{pmatrix}_{\!q}q^{(m-p)k}\gamma_{i-k}^{m-p}.$ Then, using the
definitions above, we can easily see that
$$\gamma_i^m\cdot (h\ot \gamma_j^t\vv{u,j})=\sum_{p=0}^m
\tau_{j+t}^{(p)}(\gamma_i^m)h\ot \gamma_j^{t+p}\vv{u,j}.$$

Therefore, as a left $\land$-module, $\land \ot L(u,j)$ is isomorphic
to $_{T_j(N)}(\land^N)$, that is, the left module whose action is described by $$\land \ot \land^N
\stackrel{T_j(N) \ot 1}{\rightarrow} \mat{N}{\land}\ot \land
^N\stackrel{\mbox{\tiny matrix multiplication}}{\longrightarrow}\land
^N$$ where $T_j(N): \land\rightarrow\mat{N}{\land}$ is defined by 
$T_j(N)=\begin{pmatrix} \tau_j^{(0)}&\tau_j^{(1)}& \ldots &
  \tau_j^{(N-1)}\\&\tau_{j+1}^{(0)}&&\vdots\\
&0&\ddots & \vdots\\ &&&\tau_{j+N-1}^{(0)}\end{pmatrix}.$
Note that if $\dim L(u,j)=1$, this describes a twisted module  $_\tau
(\land)$, where $\tau$ is an automorphism of $\land$.

\ 

We proceed in a similar way for the right coaction. 
Let $R_j(N)$ be an element in $\mat{N}{\land}$. Define a right
coaction on $\land^N$ as follows:
$$\land^N\stackrel{\Delta}{\longrightarrow}(\land\ot\land)^N\stackrel{(1\ot R_j(N))\cdot}{\longrightarrow}(\land\ot\land)^N\cong\land^N\ot\land.$$

Setting   $R_j(N)=
\begin{pmatrix}\rho^{(0)}_{j+N-1}\\\rho_{j+N-2}^{(1)}&\rho_{j+N-2}^{(0)}&0\\
\vdots &&\ddots\\\rho_j^{(d+N-1)}&&&\rho_j^{(0)}  \end{pmatrix}$
with $\rho_k^{(0)}=\sum_{i\in\zz_n}q^{-ik}e_i$ and $$\rho_k^{(m)}= \begin{pmatrix}-k+N+j-1-m\\m \end{pmatrix}_{\!q}\left(\prod_{p=-k+j-m}^{-k+j-1}(1-q^{p+m})\right)\displaystyle{\sum_{i\in\zz_n}
  q^{-ik}\gamma_i^m}
\mbox{ for $m\pgq 1$}$$  finally gives:

\bp The Hopf bimodule corresponding to the $\dland$-module $L(u,j)$ is
$$_{T_j(N)}^{\hspace{4pt}1}[(\land)^{\dim L(u,j)}]_1^{R_j(N)}.$$\ep

Note that when $\dim L(u,j)=1$, we have $_{T_j(1)}^{\hspace{4pt}1}[(\land)^{\dim L(u,j)}]_1^{R_j(1)},$ that is,
the usual Hopf bimodule $\land$ twisted on the left by the algebra
automorphism $T_j(1)$ and on the right by the coalgebra automorphism
$R_j(1)\cdot$ (in this case $R_j(1)=\rho_j^{(0)}$ is a grouplike element).  In particular, the trivial module $L(0,0)$ gives the
usual Hopf bimodule $\land.$

\subsection{Hopf bimodule corresponding to $P(u,j)$}

We assume here that $P(u,j)$ is not simple, that is, $N=\dim L(u,j)\neq d.$ We then have a basis for $P(u,j)$ given by $\set{A_\ell=\gamma^\ell_jD_{uj}, B_\ell=\gamma^\ell_{j+N}X^{d-N}_{uj}\,\mid\, 0\ppq \ell\ppq d-1}.$ It is easy to work out the actions of $G$ and of paths on these basis elements, and using Lemma \ref{commut}, we can also describe the action of $X$. We then obtain:

\bp The Hopf bimodule corresponding to the $\dland$-module $P(u,j)$ is
$$_{U_j}^1[(\land)^{2d}]_1^{V_j},$$ where
$U_j=\begin{pmatrix}T_j(d)&0\\0&T_{j+N}(d)  \end{pmatrix}$ and
$V_j=\begin{pmatrix} R_j(d)&0\\\Omega_j&\Pi_{j} \end{pmatrix}$ for
some matrices $\Omega_j$ and $\Pi_{j}$ which can be explicitly calculated, with $\Pi_j$ lower triangular.
 \ep


\appendix
\section{Appendix: Proof of the classification of the indecomposable
  $\dland$-modules}

We wish to describe all the indecomposable modules over
$\dland$ up to isomorphism.  We may restrict our  study to a nonsemisimple block $\Bb$ of $\dland$.
 To simplify notation, we set $m:=\frac{2n}{d}$ to be the number of
 simple modules over $\Bb,$ which we denote by $S_p:=L(u,\sigma_u^p(i))$ with $p\in
 \zz_m.$ Let $\epsilon_p$ be the corresponding idempotent, and $P_p$
 the projective cover of $S_p.$ 

We first describe the Loewy structure of the indecomposable modules,
then describe the indecomposable modules themselves. Finally, we give
a different but more elegant proof of the fact that indecomposable
modules of odd length are syzygies of simple modules (this proof is
valid for any $m$, not necessarily even).

\subsection{Loewy structure of indecomposable modules}

\bl\label{Loewy} Let $\Bb$ be a fixed block as above, with $m$ simple modules where $m\pgq 2$. Then for any
indecomposable $\Bb$-module $M$ which is not simple or projective, the
radical of $M$ is equal to the socle of $M$.
\el

\bpf 
All the indecomposable projective modules have Loewy length three, so an indecomposable non-projective module
has Loewy length at most 2, and equal to 2 if it is not simple.

Moreover $\soc(M) \subseteq \rad(M)$ in general: indeed, if not, there
would be an element $x$ in the socle which is not in the
radical. Choose $x$ such that $\Bb x$ is simple. Take a maximal
submodule $\frak{M}$ of $M$ which does not contain $x$; then
$\frak{M}\cap M=0$ so $M=\frak{M}\oplus \Bb x$ is decomposable.

Therefore,  if the Loewy length of $M$ is $2$, then  $\soc(M)$ and
$\rad(M)$ must be equal.
\epf

We now describe the top and the radical of an indecomposable module.

\bl\label{topssocles} Suppose $M$ is indecomposable and not projective or simple. One of the following must hold:

\begin{enumerate}[(i)]
\item\label{even} $M/\rad(M) = \oplus_p \epsilon_{2p}M/\rad(M)$ and $\rad(M) = \oplus_p
  \epsilon_{2p+1}\rad(M)$;
\item\label{odd} $M/\rad(M) = \oplus_p \epsilon_{2p+1}M/\rad(M)$ and $\rad(M) = \oplus_p \epsilon_{2p}\rad(M)$. 
\end{enumerate}
\el

\bpf  Choose a vector space complement $C$ with $M=\rad(M)\oplus C$ such that
$\epsilon_pC\subseteq C$ for all $p$ (this is possible since the algebra generated by the idempotents
$\epsilon_p$ is semi-simple). Then if $\alpha$ is any arrow starting at $p$ 
we have $\alpha C = \alpha \epsilon_pC \subseteq \epsilon_{p+1}\rad(M) \oplus \epsilon_{p-1}\rad(M)$. But $p\pm 1$ have the same
parity and $m$ is even, so if
$$M' := [\oplus_p \epsilon_{2p}C ] \oplus [\oplus _p \epsilon_{2p+1}\rad(M)], \ \ 
M'' := [\oplus_p \epsilon_{2p+1}C]\oplus [\oplus_p \epsilon_{2p}\rad(M)]
$$
then these are $\Bb$-submodules of $M$ (use Lemma \ref{Loewy}) and $M=M'\oplus M''$. 
If $M=M'$ then we have (\ref{even}) and otherwise (\ref{odd}).
\epf

\bd We say that the top of $M$ is {\bf even} if Lemma \ref{topssocles}(\ref{even}) holds, and
that the top of $M$ is {\bf odd} if Lemma \ref{topssocles}(\ref{odd}) holds.

\ed

\subsection{Description of the indecomposable modules over a block $\Bb$}

\subsubsection{Applying Fitting's Lemma}

Let $M$ be an indecomposable module.
We can assume, without loss of generality, that the top of $M$ is
even. Set $M_p=\epsilon_pM:$ this is $\epsilon_pM/\rad(M)$ if $p$ is even, and $\epsilon_p \rad(M)$ otherwise.

 Then $M$ is completely
described by the linear maps given by left multiplication,
$$b_{2r}: M_{2r} \to M_{2r+1}, \ \  \overline{b}_{2r-1}: M_{2r}\to M_{2r-1}
$$ (the other arrows act as zero).
We would like to have maps which we can compose, so  we define $\beta_p:=\, {^t}(\overline{b}_p): M_p^*\to M_{p+1}^*$,
for $p$ odd. For any fixed basis of $M_p$ and $M_{p+1}$, with respect to the dual bases,
the matrix of the map $\beta_p$ is the transpose of the matrix of
$\overline{b}_p$.  Let $b$ and $\beta$ also denote the matrices of the maps
$b$ and $\beta$.

Consider the map $\alpha_p: M_p \to M_p$, which is equal to
$\beta_{p-1}b_{p-2}\ldots b_{p+2}\beta_{p+1}b_p$ if $p$ is even, and to
$b_{p-1}\beta_{p-2}\ldots \beta_{p+2}b_{p+1}\beta_p$ if $p$ is odd, that is,   the
composition of all $b$'s and $\beta$'s,  one each 
in the only possible order.

By Fitting's Lemma, there is some large $N>0$ such that 
$U_p:= \ker(\alpha_p^N) = \ker(\alpha_p^{N+\ell})$ and $V_p:= \im(\alpha_p^N)=\im(\alpha_p^{N+\ell})$
for all $\ell\geq 0$, where $U_p$ and $V_p$ are $\alpha_p$-invariant and $M_p = U_p \oplus V_p$. Moreover $\alpha_p|_{U_p}$ is nilpotent and $\alpha_p|_{V_p}$ is an
isomorphism.

Set $U:= \oplus _p U_p$ and $V:= \oplus_p V_p$. Then $M=U\oplus V$ and moreover, both
$U$ and $V$ are invariant under the actions of  $b_{2r}$ and
$\beta_{2r+1}$. Let us check this for $b_{2r}:$ note that
$\alpha_{2r+1}b_{2r}=b_{2r}\alpha_{2r}.$ So if $x$ is in $U_{2r}$,
then $\alpha_{2r+1}^Nb_{2r}(x)=b_{2r}\alpha_{2r}^N(x)=0,$ so
$b_{2r}(x)$ is in $U_{2r+1}.$ If $x$ is in $V_{2r},$ then we can write
$x=\alpha_{2r}^N(z)$ and hence $b_{2r}(x)=\alpha_{2r+1}^N(z)\in
V_{2r+1}.$ 

The decomposition $M=U\oplus V$ above is a decomposition of
$\Bb$-modules: for each $M_p$ take a basis for $U_p$ and a basis for
$V_p$ such that the union is a basis for $M_p$. Then, for each $t$ odd,
the matrix of $\beta_t$  is a block-diagonal matrix. Now return to
$\bar{b}_t$ (instead of $\beta_t$): we take the transpose of this
matrix, which has the same block form. Therefore $M=U\oplus V$ as a
$\Bb$-module.

\subsubsection{Canonical forms}\label{canonicalforms}

We now assume that $M$ is an indecomposable module. Therefore, $M$ is
equal to $U$ or $V.$

\paragraph{First case: $M=U$.} Consider the algebra $A$ generated by
the $\beta_{2p-1}$ and the $b_{2p}.$ This is a Nakayama algebra, whose
quiver is cyclic, and all the paths of length at least $\delta$ where
$\delta=m\cdot{\rm max}\set{\mbox{nilpotence degree of the }\alpha_p}$
are zero, with $\delta\pgq \dim M.$

We know that modules over $A$ are uniserial. Therefore there is a
generator $v$ of $M$ over $A$ such that $v=\epsilon_p v,$ and a
basis for $M$ is given by $\set{v,b_pv,\beta_{p+1}b_pv,\ldots}$
(the index $p$ is even here, since we are in the even top case). 

Returning to $b$'s
 and $\bar{b}$'s, we can write the matrices for the actions of these
 arrows, and we see that $M$ is a string module. If it has odd length
 then we recover a syzygy of a simple module, and otherwise we have a
 string module of type $M_2^+(u,i)$. Note that if we had started with
 a module with odd top, we would have obtained either a syzygy of a
 simple module (if the length were odd), or  a string
 module of type $M_2^-(u,i)$ (if the length were even). See Section
 \ref{desc} for the notations.

\paragraph{Second case: $M=V$.} We have $\ker(\alpha_p^N)=0$ for all
$p$ and $M=\bigoplus_p\im(\alpha_p^n)$, so we deduce that the maps
$b_{2r}$ and the maps $\beta_{2r-1}$ (and hence $\bar{b}_{2r-1}$) are
invertible. Let $\ell$ be the dimension of the $M_p$ (they all have
the same dimension). We can choose bases for the $M_p$ such that all
the maps $b_{2r}$ and $\bar{b}_{2r-1}$ are represented by the identity
matrix, except for one of them, say $b_0.$ Note that the bases for the
$M_p$ are determined by a choice of basis for $M_1$ say, and any
change of basis for $M_1$ leads to the same change of basis for all
the other $M_p.$ Since $M$ is
indecomposable, the map $b_0$ must be an indecomposable linear map, so
we can choose a basis for $M_1$ such that the matrix of $b_0$ is of
the form $J_\ell(\lambda)=\mx{\lambda & 1
  &\ldots\\0&\lambda&1&\ldots\\&&\cdots\\0&\ldots&&\lambda}$ for some
nonzero $\lambda$ in $k,$ and such that all the other matrices remain
identity matrices. 

We call this module $C^{\ell +}_\lambda(u,i).$ The same construction
starting with a module whose top is odd will give us $C^{\ell
  -}_\lambda(u,i).$ Note that they have even length, equal to
$m\ell.$

\subsubsection{Auslander-Reiten components}

We refer the reader to  \cite[V.1 and VII.1]{ARS} and \cite[I.7 and I.8]{erdmann} for information about Auslander-Reiten sequences and quivers.

\begin{enumerate}[{\bf (I)}]

\item {\bf Indecomposable modules of odd length:}  For each simple module $S_p$ we have the standard Auslander-Reiten sequence
$$0 \to \Omega(S_p) \to {\rm rad}(P_p)/{\rm soc}(P_p) \oplus P_p \to \Omega^{-1}(S_p)\to 0
$$
These are the only Auslander-Reiten sequences where projectives occur.
For a  symmetric algebra (which is the case here, see  the Introduction), the Auslander-Reiten translation is the same as $\Omega^2$, so it
follows from Section \ref{canonicalforms} that any indecomposable module of odd length, say $M = \Omega^k(S_p)$
with $k\neq -1,$ has Auslander-Reiten sequence
$$0 \to \Omega^{k+2}(S_p) \to \Omega^{k+1}(S_{p-1})\oplus \Omega^{k+1}(S_{p+1}) \to M \to 0,
$$ obtained by applying $\Omega^{k+1}$ to the sequence above.

This describes all the  Auslander-Reiten sequences, and therefore the
Auslander-Reiten component in which $M$ occurs.

\item {\bf String modules of even length: } By explicit calculation,
  we can easily see that
  $\Omega(M_{2\ell}^{+}(u,i))=M_{2\ell}^+(u,\sigma_u^{-1}(i)).$ So
  $M_{2\ell}^{+}(u,i)$ has $\Omega$-period equal to $m,$ hence it
  lies in a tube of rank $\frac{m}{2}=\frac{n}{d}.$

We now determine the Auslander-Reiten sequences. Let us start with the
Auslander-Reiten sequence for $M^+_2(u,i):$ it is of the form
$0\rightarrow M^+_2(u,i-d)\rightarrow M \rightarrow
M^+_2(u,i)\rightarrow 0$ since the Auslander-Reiten translate is given
by $\Omega^2.$ The module $M$ must have length four, and cannot be a
direct sum of modules of length two (otherwise the exact sequence
would split). The composition factors of $M$ are
$L(u,\sigma_u^{-2}(i)),$ $L(u,\sigma_u^{-1}(i)),$ $L(u,i)$ and
$L(u,\sigma_u(i)).$ Moreover, $L(u,\sigma_u^{-2}(i))$ is in the socle
of $ M^+_2(u,i-d)$ so it must be in the socle of $M$, and $L(u,i)$ is
in the top of $M^+_2(u,i)$ so it must be in the top of $M.$ It then
follows from the classification of the modules in Section
\ref{canonicalforms} that $M$ must be $M_4^+(u,i-d).$

Let us now determine the Auslander-Reiten sequence for  $M_4^+(u,i).$
We know that the left hand term of the sequence must be
$\Omega^2(M_4^+(u,i)),$ and we know exactly which modules of length two
occur, from the Auslander-Reiten sequences for modules of length two. Therefore we
have $0\rightarrow M_4^+(u,i-d)\rightarrow M_2^+(u,i)\oplus N
\rightarrow M_4^+(u,i) \rightarrow 0$, with $N$ indecomposable of length six. The
composition factors of $N$ are $L(u,\sigma_u^{-1}(i)),$
$L(u,\sigma_u(i))$, $L(u,\sigma_u^{3}(i))$, $L(u,\sigma_u^{-2}(i))$
$L(u,i)$ $L(u,\sigma_u^{2}(i)).$ The first one must be in the socle,
and hence so must the next two. The last one must be in the top, and
therefore so are the two remaining simple modules. Therefore
$N=M_6^+(u,i).$

The general sequences in Section \ref{arcomponents}{\bf
  (\ref{arstring})} are obtained in this way by induction. Therefore,    the $M_{2\ell}^{+}(u,\sigma_u^p(i))$ for all
even $p$ and all $\ell$ form precisely one tube, where the modules of
length $2\ell$ form the $\ell^{th}$ row. The $M_{2\ell}^{+}(u,\sigma_u^p(i))$ for all
odd $p$ and all $\ell$ form precisely one tube, where the modules of
length $2\ell$ form the $\ell^{th}$ row.

Similarly, there are two tubes, each of rank $\frac{n}{d}, $
containing all the modules $M_{2\ell}^{-}(u,\sigma_u^p(i))$.

\item {\bf Band modules:} We can see that $\Omega(C^{\ell
    +}_\lambda(u,i))\cong C^{\ell -}_{-\lambda}(u,i)$ and $\Omega(C^{\ell
    -}_\mu(u,i))\cong C^{\ell +}_{-\mu}(u,i).$ Hence $C^{\ell
    +}_\lambda(u,i)$ has period two and lies in a tube of rank one. 

We indicate briefly how to determine the Auslander-Reiten sequences
given in Section \ref{arcomponents}{\bf (\ref{arband})}: the
Auslander-Reiten sequence for $C_\lambda^{1+}(u,i)$ is of the form $0
\rightarrow C_\lambda^{1+}(u,i) \rightarrow M \stackrel{\pi}{\rightarrow}
C_\lambda^{1+}(u,i) \rightarrow 0$ (recall that  the Auslander-Reiten
translate is given by $\Omega^2$ and that these modules are periodic
of period two). The module $M$ must be indecomposable of length $\frac{4n}{d}$, its top must
be $\bigoplus_pL(u,\sigma_u^{2p}(u,i))^2$ and its socle is
$\bigoplus_pL(u,\sigma_u^{2p+1}(u,i))^2$ (as vector spaces), so $N$ is
of the form $C_\mu^{2+}(u,i).$ We need to determine $\mu.$ Let
$\set{v_1,v_2}$ be a basis for $L(u,i)^2$ such that the action of
$b_0$ is given by the matrix $\mx{\mu&1\\0&\mu}$. The vector $v_1$
generates a submodule so must come from the left hand term in the
sequence. Hence $\pi(v_1)=0$, and $\pi(v_2)\neq 0$ is a basis for
$L(u,i)$ in $C_\mu^{2,+}(u,i)$. We then have
$\pi(b_0v_2)=\pi(v_1+\mu v_2)=\mu \pi(v_2),$ and
$\pi(b_0v_2)=b_0\pi(v_2)=\lambda\pi(v_2)$ since the action of $b_0$ on
$L(u,i)$ is given by multiplication by $\lambda$. Hence $\mu=\lambda.$
An induction gives the expected sequences, and therefore, 
 for each $\lambda$, the modules $C^{\ell
  +}_\lambda(u,i)$ form one tube, with $C^{\ell +}_\lambda(u,i)$ in
row $\ell.$ Similarly, the $C^{\ell -}_\lambda(u,i)$ when $\ell$
varies form one tube,
with $C^{\ell -}_\lambda(u,i)$ in row $\ell.$

\end{enumerate}

\subsection{Indecomposable modules of odd length: another proof}

We work over a block $\Bb$ with $m$ simple modules, where $m$ is any
integer at least equal to 2.

\bl\label{oddlength}  Suppose $M$ is indecomposable in $\Bb$ of odd length. Then some syzygy $\Omega^k(M)$
for some $k \in \mathbb{Z}$ is simple. \el

\bpf Since $M$ has odd length, $\soc(M)$ and $M/\soc(M) (= \Top(M))$ have different length, say 
$\soc(M)$ has larger length. Consider a projective cover
$$0 \to \Omega(M) \to P \to M \to 0.$$
We know from the structure of projective modules that $\Top(P)
\cong \soc(P)$, so  $\Top(M) \cong \soc(\Omega(M))$.

If $\Top(M)$ has length $t$, then
$P$ has length $4t$ (since $\Top( P) =\Top( M)$). It follows that 
$\soc(\Omega(M))$ has length $t$ and $\Top(\Omega(M))$ has length $t+x$
if  $\Omega(M)$ is not simple.

Therefore the length of the module $\Omega(M)$ is strictly smaller
than that of $M,$ and $\Top (\Omega(M))$ has more composition factors
than $\soc(\Omega(M))$. If $\Omega(M)$ is simple then the proof is
finished. Otherwise we can repeat the argument, and since the lengths
are positive integers, the process stops  after a finite
number of steps, that is, some $\Omega^k(M)$ is simple.

If the length of $\soc (M)$ is smaller than the length of $\Top(M),$ a
similar argument using injective envelopes gives the result.
\epf

\br It follows from the proof above that $\abs{{\rm length}(\Top(M)) - {\rm
  length}(\soc(M))}=1$ for any indecomposable module $M$ of odd length.
\er


\flushleft{{\hrulefill\hspace*{10cm}}\\
\scriptsize{ 
{\sc Karin Erdmann\\Mathematical Institute,\\
24-29 St. Giles,\\
Oxford OX1 3LB,\\
United Kingdom.\\ E-mail:} erdmann@maths.ox.ac.uk\\
\vspace*{.3cm}
{\sc Edward L. Green\\Department of Mathematics,\\
Virginia Polytechnic Institute and State University,
\\Blacksburg, VA 24061-0123\\ USA.\\ E-mail:} green@math.vt.edu\\ 
\vspace*{.3cm}
{\sc Nicole Snashall\\Department of Mathematics,\\ University of Leicester,\\ Leicester LE1 7RH, \\United
    Kingdom.\\ E-mail:} N.Snashall@mcs.le.ac.uk\\ 
\vspace*{.3cm}
{\sc Rachel Taillefer\\Laboratoire d'Arithm\'etique et d'Alg\`ebre,\\Facult\'e des Sciences et Techniques,\\23 Rue Docteur Paul Michelon,\\42023 Saint-Etienne Cedex 2,\\France.\\ E-mail:} rachel.taillefer@univ-st-etienne.fr\\{\sc Telephone:+33 (0)4 77 48 15 33}\\{\sc Fax:+33 (0)4 77 48 51 33.}

}

\end{document}